\documentclass[11pt,a4paper]{article}
\usepackage{bbm}
\usepackage{mathrsfs}
\usepackage{amsfonts}
\usepackage{}

\usepackage{graphicx}
\usepackage{amsfonts,amsmath, amssymb}
\usepackage{color}

\newtheorem{theo}{\textbf{\ \ \quad Theorem}}[section]

\newtheorem{lem}{\textbf{\ \ \quad Lemma}}[section]
\newtheorem{remark}{\textbf{\ \ \quad Remark}}[section]

\newtheorem{prop}{\textbf{\ \ \quad Proposition}}[section]
\newtheorem{defi}{\textbf{\ \ \quad Definition}}[section]

\newcommand{\lbl}[1]{\label{#1}}

\newcommand{\be}{\begin{equation}}
\newcommand{\ee}{\end{equation}}
\newcommand\bes{\begin{eqnarray}}
\newcommand\ees{\end{eqnarray}}
\newcommand{\bess}{\begin{eqnarray*}}
\newcommand{\eess}{\end{eqnarray*}}

{\setlength\arraycolsep{2pt}}

\title{Periodic solution of stochastic process in the distributional sense}

\author{Guangying Lv$^a$, Hongjun Gao$^c$ and Jinlong Wei$^d$\\
\\
\ \\
   {\small \it $^a$Institute of Applied Mathematics, Henan University}\\
  {\small \it Kaifeng, Henan 475001, China}\\
  {\small\it $^b$Center for Applied Mathematics, Tianjin University}\\
  {\small \it Tianjin 300072, China}\\
  {\small \tt gylvmaths@henu.edu.cn}\\
    {\small \it $^c$Institute of Mathematics, School of Mathematical Science}\\
  {\small \it Nanjing Normal University, Nanjing 210023, China}\\
  {\small \tt gaohj@njnu.edu.cn}\\
  {\small \it $^d$ School of Statistics and Mathematics, Zhongnan University of}\\
  {\small \it
 Economics and Law, Wuhan, Hubei 430073, China}\\
   {\small \tt  weijinlong.hust@gmail.com }
}

\begin{document}
\maketitle

\medskip

\begin{abstract} In this paper, we aim to study a stochastic process from a macro point of view, and
thus periodic solution of a stochastic process in distributional sense is
introduced. We first give the definition and then establish
the existence of periodic solution on bounded domain.
Lastly, for the case that probability density function
exists, we obtain the existence periodic solutions of the probability
density function corresponding to the stochastic process
by using the technique of deterministic partial differential equations.



\end{abstract}

\baselineskip=15pt

\section{Introduction}
\setcounter{equation}{0}

Some properties of a stochastic process are worth being studied, such as
the long time behavior, periodicity, ergodicity and so on. There are many classical
theories, see the books \cite{IW,O,Schuss}.
In this paper, we will give a new viewpoint about the periodicity
of stochastic processes. We consider the stochastic process from another fact--a macro point of
view. We do not consider the motion of a single particle, while we are
concerned with the motion the entire system. It is well-known that
the probability density function (PDF) of a stochastic process can describe
the entire distribution of a system. Hence,
in this paper, we consider some property of entire system--time-periodicity of PDF. The density of a stochastic
process is called as Fokker-Planker equation or Kolmogorov equation, which has been studied in
\cite{BKRS}.

Now, we consider the multidimensional Fokker-Planck equation for the following SDE
  \bes
dX_t=b(t,X_t)dt+\sigma(t,X_t)dB_t, \ \ \ X_0=x\in\mathbb{R}^d,
  \lbl{1.1}\ees
where $b$ is an $d$-dimensional vector function, $\sigma$ is an $d\times m$ matrix function and
$B_t$ is an $m$-dimensional Brownian motion, see Page 99 in \cite{Duan2015}.
Let the probability density function of (\ref{1.1}) be $p(t,x)$ if it exists, then we deduce that
$p(t,x)$ satisfies
   \bes
\partial_tp(t,x)=\frac{1}{2}div(div(\sigma\sigma^Tp))-div(b(t,x)p)
  \lbl{1.2} \ees
with initial data
  \bes
p(0,x)=p_0(x), \ \ x\in\mathbb{R}^d.
  \lbl{1.3}\ees
In this paper we mainly consider the property of $p(t,x)$. There is a fact that
for most of stochastic process, the PDF may not exist. Therefore, in this paper,
we first give the notion of periodic solution in distributional sense for discrete time and
continuous time stochastic process and then consider
a special case, in which the PDF exists.

Let us recall some development of periodic solutions.
Periodic solutions have been a central concept in the theory of the deterministic dynamical system
starting from Poincar\'{e}'s work \cite{Poincare1988}. For a random periodic dynamic system, to study
the pathwise random periodic solutions is of great importance. Zhao-Zheng \cite{Zhao2009} started to study
the problem and gave a definition of pathwise random periodic solutions for $C^1$-cocycles. Recently, Feng et al. did
some beautiful work about the periodic solutions, see \cite{Feng2011jde,Feng2012jfa,Feng2016jfa}.
Noting that the definition of periodic solution in \cite{Feng2011jde,Feng2012jfa,Feng2016jfa,Zhao2009} is different
from here, we consider the time-periodicity of entire system in distributional sense.
During we prepare our paper, we find the paper of
Chen et al. \cite{Chenfeng2017}, where the existence of periodic solutions of Fokker-Planck equations is considered.
They obtained the desired results by discussing the existence of periodic solutions in distributional sense for some
stochastic differential equations (SDEs). More precisely, they used the properties of solutions of SDEs to study the
properties of solutions of Fokker-Planck equation. They obtained the time-periodicity of PDF in the whole space.
We will give another proof in viewpoint of PDEs. Moreover, the definition of
periodic solution to discrete time and discrete state stochastic process will be given.
The topic of periodic solutions to stochastic process in distributional sense on bounded domain
is also considered in this paper. About the almost periodic solution, see
\cite{LS2014jfa,WL2012}.

The rest of this paper is arranged as follows. In Sections 2, we present some known results on
PDEs' theory.  In Sections 3, we first give some definitions of periodic solutions
to stochastic process in distributional sense, then establish the
existence of periodic solution on bounded domain by using the
method of \cite{Chenfeng2017}. For the case that the PDF exists,
we obtain the existence of periodic solution of Fokker-Planck equations
on bounded domain and in the whole space by using the
method of deterministic partial differential equations in Section 4.

\section{Some known results}\label{sec2}\setcounter{equation}{0}
In this section, we recall some known results about existence of PDF of the diffusion
It\^{o} process and the existence of periodic solution of parabolic equations.

 Consider a Markov process in $\mathbb{R}^d$ with transition probabilities $P(s,x,t,B)$ ($B$ is a
Borel set in $\mathbb{R}^d$)
is called a \emph{diffusion process} or a \emph{diffusion} if there is a mapping $b:\mathbb{R}^d\times[0,\infty)\to\mathbb{R}^d$,
called the \emph{drift coefficient}, and a mapping $(x,t)\mapsto A(x,t)$ with values in the space of symmetric
operator on $\mathbb{R}^d$, called the \emph{diffusion coefficient} or \emph{diffusion matrix}, such that

(i) for all $\varepsilon>0$, $t\geq0$ and $x\in\mathbb{R}^d$ we have
   \bess
\lim\limits_{h\to0}h^{-1}P(t,x,t+h,V(x,\varepsilon))=0,
   \eess

(ii) for some $\varepsilon>0$ and all $t\geq0$, $x\in\mathbb{R}^d$ we have
     \bess
\lim\limits_{h\to0}h^{-1}\int_{U(x,\varepsilon)}(y-x)P(t,x,t+h,dy)=b(x,t),
   \eess

(iii) for some $\varepsilon>0$ and all $t\geq0$, $x,z\in\mathbb{R}^d$ we have
     \bess
\lim\limits_{h\to0}h^{-1}\int_{U(x,\varepsilon)}\langle y-x,z\rangle P(t,x,t+h,dy)=2\langle A(x,t)z,z\rangle,
   \eess
where $U(x,\varepsilon)=\{y:|x-y|<\varepsilon$\} and $V(x,\varepsilon)=\{y:|x-y|>\varepsilon$\}.
If $A$ and $b$ do not depend on $t$, then the diffusion is homogeneous. Bogachev et al. \cite{BKRS} obtained
the following proposition.
\begin{prop}\lbl{p1.1} Suppose that relations (i)-(iii) hold locally uniformly in $x$ and the
functions $a^{ij},b^i$ ($A=(a^{ij}),\,b=(b_j)$) are locally bounded. Then the transition probabilities satisfy the parabolic
Fokker-Planck-Kolmogorov equation
   \bess
\partial_t\mu=\partial_{x_i}\partial_{x_j}(a^{ij}\mu)-\partial_{x_i}(b^i\mu)
   \eess
in the sense of generalized functions. If $\nu$ is a finite Borel measure on $\mathbb{R}^d$
and
   \bess
\mu_t(dx)=\int_{\mathbb{R}^d}P(0,y,t,dx)\nu(dy),
   \eess
then the measure $\mu=\mu_t(dx)dt$ gives a solution to the Cauchy problem with the
initial condition $\nu|_{t=0}=\nu$.
  \end{prop}
The above proposition is concerned with the case of measure-valued solution about the
Fokker-Planck equation. The next result shows that there exists PDF for a stochastic
process under some assumptions.
Let $D_T=\Omega\times(0,T)$, where $D\subset\mathbb{R}^d$ is an open
set and $T>0$ is a fixed number. Bogechev et al. \cite{BKRS} obtained the following result.
\begin{prop}\cite[Theorem 6.3.1]{BKRS} Let $\mu$ be a locally finite Borel measure on $D_T$ such that
$a^{ij}\in L^1_{loc}(D_T,\mu)$ and
   \bess
\int_{D_T}\left[\partial_t\phi+a^{ij}\partial_{x_i}\partial_{x_j}\phi\right]d\mu\leq C(\sup_{D_T}|\phi|+\sup_{D_T}|\nabla_x\phi|)
   \eess
for all nonnegative $\phi\in C_0^\infty(\Omega_T)$. Then the following assertions are true.

(i) If $\mu\geq0$, then $(det A)^{1/(d+1)}\mu=\rho dxdt$, where $\rho\in L_{loc}^{(d+1)'}(D_T)$.

(ii) If, on every compact set in $D_T$, the mapping $A$ is uniformly bounded, uniformly
nondegenerate, and H$\ddot{o}$lder continuous in $x$ uniformly with respect to $t$, then
$\mu=\rho dxdt$, where $\rho\in L_{loc}^r(D_T)$ for every $r\in[1,(d+2)')$.
\end{prop}

The above Proposition is the existence of probability density function in the whole space. Now, we consider the
bounded domain. As stated in \cite{Duan2015}, in the simulations, we have to take $x$ in a large but bounded domain
$D\subset\mathbb{R}^d $ and we could impose absorbing boundary condition on $\partial D$, i.e., as long as
a "particle" or a solution path reaches the boundary, it is removed from the system. The above assumptions
implies the following system
   \bes\left\{\begin{array}{llll}
  \partial_tp(t,x)=A^*p(t,x), \ \ t>0,\ x\in D,\\
p|_{\partial D}=0,\\
p(0,x)=p_0(x), \ \ \ x\in D,
   \end{array}\right.\lbl{2.1}\ees
where
   \bess
A^*p=-\sum_i\frac{\partial}{\partial x_i}(b_ip)+\frac{1}{2}\sum_{i,j}
\frac{\partial^2}{\partial x_i\partial x_j}((\sigma\sigma^T)_{ij}p).
  \eess
Due to the absorbing boundary condition, the particle will not come back when it reach the boundary.
Thus under the absorbing boundary, it is impossible to get the existence of periodic solution to (\ref{2.1}).
Therefore, we must consider another  case: the reflecting boundary condition \cite[Section 5.1.1]{Gardiner2009book}.
The Fokker-Planck equation of (\ref{1.1}) can be written as
   \bes
  \partial_tp(t,x)=\nabla\cdot\left[\frac{1}{2}(\nabla\cdot(\sigma\sigma^Tp))-bp(t,x)\right].
    \lbl{2.2}\ees
The reflecting boundary condition means particles or solution paths can not
leave a bounded domain $D$, and hence there is zero net flow of $p$ crossing the boundary
$\partial D$. Thus we impose the following reflecting boundary condition
   \bes
\frac{1}{2}(\nabla\cdot(\sigma\sigma^Tp))-bp(t,x)=0 \ \ \ {\rm on}\ \ \ \partial D.
   \lbl{2.3}\ees
Integrating (\ref{2.2}) over $D$ and using the boundary condition (\ref{2.3}) together with the divergence theorem,
we have conservation of probability
   \bess
\frac{\partial}{\partial t}\int_Dp(t,x)=0,
  \eess
that is to say,
   \bess
\int_Dp(t,x)=\int_Dp_0(x)=1.
   \eess
In this case, it is possible to obtain the existence of periodic solution to (\ref{2.2})-(\ref{2.3}) with
initial data. In order to obtain the desired results, we recall some results about the periodic parabolic
equations, see \cite{Hess1991book}.

Now, we consider the periodic-parabolic eigenvalue problem
   \bes\left\{\begin{array}{llll}
   \partial_tu+\mathcal {A}(t)u=\mu u\ \ \ & in\ \ D\times\mathbb{R},\\
   \mathcal {B}u=0\ \ \ \ & on \ \ \partial D\times\mathbb{R},\\
   u \ \ T-periodic \ in\ t,
\end{array}\right.
 \lbl{2.4}\ees
where $\mathcal {A}(t)$ is a uniformly elliptic differential operator of second order
depending $T$-periodically on $t$, i.e.,
   \bess
\mathcal {A}(t)u=\mathcal {A}(t,x,D)u=-\sum_{j,k=1}^da_{jk}(t,x)\frac{\partial^2}{\partial{x_j}\partial{x_k}}u
+\sum_{j=1}^da_{j}(t,x)\frac{\partial}{\partial{x_j}}u+a_0(t,x)u,
   \eess
 and
   \bess
\mathcal {B}u=\left\{\begin{array}{lll}u \qquad \qquad \qquad Dirichlet \ b.c.,\\
\frac{\partial u}{\partial\nu}+b_0(x)u \ \ \ Neumann\ or \ regular\ oblique\ derivative\ b.c..
   \end{array}\right.\eess

We say $\mu\in \mathcal {C}$ ($\mathcal {C}$ denotes complex value) is an eigenvalue if there is a nontrivial
solution $u$ (eigenfunction) of (\ref{2.4}). We search in particular for an eigenvalue $\mu\in\mathbb{R}$ having
a positive eigenfunction ("principal eigenvalue" $\mu$).

In order to establish the existence of solutions of (\ref{2.4}), we consider the inhomogeneous linear evolution equation
    \bes\left\{
\begin{array}{lll}
\dot{u}(t)+\mathcal {A}(t)u=f(t), \ \ \ 0<t<T,\\
u(0)=u_0, \ \ \ \ \ u_0\in X,
   \end{array}\right.\lbl{2.5}\ees
where $f\in C^\theta([0,T],X)$, $0<\theta\leq1$ and $X$ is a Banach space.
Assume the closed linear operator $\mathcal{A}$ in $X$ satisfies

  (i) $dom(\mathcal{A}):=dom(\mathcal{A}(t))$ is dense in $X$ and independent of $t$,

  (ii) $\{\lambda\in \mathcal {C}:Re\lambda\leq0\}\subset\rho(\mathcal{A}(t))$, $\forall t\in[0,T]$, ($\rho(\mathcal{A}(t))$ denotes
  the resolvent set of operator $\mathcal{A}(t)$),

  (iii) $\|(\mathcal{A}(t)-\lambda)^{-1}\|\leq\frac{c}{1+|\lambda|}$, $\forall \lambda\in \mathcal {C}$, $Re\lambda\leq0$, $\forall t\in[0,T]$. \\
Set $\mathcal {A}:=\mathcal{A}(0)$ and take the fractional power spaces $X_\alpha$ with respect to
$\mathcal{A}$. Assume further

(iv) $\mathcal{A}(\cdot):[0,T]\to \mathcal{L}(X_1,X)$ is H\"{o}lder continuous.\\
It follows from the results of Sobolevskii \cite{Sobolevskii1966} that there exists a unique solution $u$ of (\ref{2.5}) with
   \bess
u\in C([0,T],X)\cap C^1((0,T],X)\ \ if \ \ u_0\in X.
   \eess
Moreover, there exists the evolution operator
   \bess
U(t,s)\in\mathcal {L}(X)
   \eess
such that the solution of (\ref{2.5}) can be represented in the following form
   \bes
u(t)=U(t,0)u_0+\int_0^tU(t,s)f(s)ds,\ \ \ \ 0\leq t\leq T.
   \lbl{2.6}\ees
The function $U$ is strongly continuous on the set $\triangle:=\{(t,s)\in[0,T]\times[0,T]:0\leq s\leq t\leq T\}$,
i.e. $U(\cdot)u_0\in C(\triangle,X)$ for each $u_0\in X$, and satisfies
   \bess
U(t,t)=I, \ \ \ U(s,t)U(t,\tau)=U(s,\tau), \ \ 0\leq\tau\leq t\leq s\leq T.
  \eess
Set $K:=U(T,0)$. The Krein-Rutman theorem implies that $r:=spr(K)>0$, where $spr(K)$ is the principal
eigenvalue of $K$.

Let
   \bess
\mathcal {L}:=\partial_t+\mathcal {A}(t),
   \eess
and set
    \bess
L&:=& \ the \ operator\ in \ \mathbb{F}_1\ introduced\ by \ \mathcal {L},\ \mathcal {B}\\
 &&  and \ the \ T-periodicity,\ with\ domain\ dom{L}=\mathbb{F}_1,
    \eess
where
   \bess
\mathbb{F}_1:=\{w\in C^{2+\theta,1+\frac{\theta}{2}}(\bar D\times\mathbb{R}):\ \mathcal {B}w=0\ \ on \ \partial D\times\mathbb{R},
w\ T-periodic \ in \ t\}.
   \eess
Assume the following conditions hold:

(A) $\mathcal {A}(t)$ is uniformly elliptic for each $t\in\mathbb{R}$ and $T$-periodic in $t$, of given period $T>0$. More precisely,
assume the coefficient functions $a_{jk}=a_{kj}$, $a_j,\,a_0$ belong to the space
    \bess
\mathbb{F}:=\{w\in C^{\theta,\frac{\theta}{2}}(\bar D\times\mathbb{R}): w\ T-periodic \ in \ t\}.
   \eess
We keep $\mathcal{B}=\mathcal{B}(x,D)$ independent of $t\in[0,T]$, such that the operator $\mathcal{A}(t)$, the realization of
$(\mathcal{A}(t),\mathcal{B})$ in $L^p(D)$ ($N<p<\infty$) has domain independent of $t$. We
assume that
    \bess
a_0(t,x)\geq1, \ \ \ \ \forall (t,x)\in[0,T]\times \bar D.
   \eess
 Then $\{\mathcal{A}(t):0\leq t\leq T\}$ satisfies the hypotheses (i)-(iv).  Thus, by
 the results of Sobolevskii \cite{Sobolevskii1966}, we get the existence of evolution operator
 $U(t,s)$ for $0\leq s\leq t\leq T$.

Now, we give the relation between the solutions of (\ref{2.4}) and (\ref{2.5}) with $f=0$.

\begin{prop}\lbl{p2.3} Assume $(A)$ holds. Then we have :

$r:=spr(K)$ is principal eigenvalue of $K$, with principal eigenfunction $u_0\gg0$$\Leftrightarrow$
$\mu=-\frac{1}{T}\log r$ is an eigenvalue of $L$ with positive eigenfunction $u:=u(t)=e^{\mu t}U(t,0)u_0$.
  \end{prop}

We have the following proposition about the positivity of $\mu$.
\begin{prop}\lbl{p2.4} Assume $(A)$ holds. Assume further that
the zero-order term of $\mathcal {A}(t)$ satisfies $a_0\geq0$ on $\bar D\times\mathbb{R}$,
and that
   \bess
a_0\not\equiv0\ \ on\ \bar D\times\mathbb{R}\ \ \  \ if \ \mathcal {B}=\frac{\partial}{\partial\nu}.
   \eess
Then $0<r<1$.
  \end{prop}

\section{Definitions of periodic solutions in distributional sense}\label{sec3}\setcounter{equation}{0}
In this section, we give some definitions of periodic solutions in distributional sense, including
discrete time and discrete state stochastic process (also called stochastic sequence)
and continuous time and continuous state stochastic process.

We start to consider the discrete time and discrete state stochastic process. Suppose
a stochastic sequence $\{X_n\}_{n\geq1}$ defined on a complete probability space has a
one-step transition probability matrix $P$. Following the Chapman-Kolmogorov equation,
we have the $N$-th step transition probability matrix $P^{(N)}$ satisfying
   \bess
P^{(N)}=P\cdot P^{(N-1)}=\cdots=P^N.
   \eess
Now, we suppose each particle has $m$ state
 in a particle system and the particle system has an initial distribution
$(x_1^0,x_2^0,\cdots,x_m^0)^T$. Consider the distribution of the system after being
transferred $N$-step (denoted by $(x_1^N,x_2^N,\cdots,x_m^N)^T$),  we have
   \bess
(x_1^N,x_2^N,\cdots,x_m^N)^T=P^{(N)}(x_1^0,x_2^0,\cdots,x_m^0)^T=P^N(x_1^0,x_2^0,\cdots,x_m^0)^T.
   \eess
Therefore, if the following holds
   \bess
P^{(N)}(x_1^0,x_2^0,\cdots,x_m^0)^T=(x_1^0,x_2^0,\cdots,x_m^0)^T,
  \eess
then the particle system turn back to the initial distribution.
We give the first definition of periodic
solution in distributional sense.

\begin{defi}\lbl{d4.1} (discrete time and discrete state stochastic process) Suppose a particle
system has one-step transition probability matrix $P$ and contains $m$ states with the  initial distribution
$(x_1^0,x_2^0,\cdots,x_m^0)^T$. If there exists a positive constant
$N\in\mathbb{N}$ such that
   \bes
P^{(N)}(x_1^0,x_2^0,\cdots,x_m^0)^T=(x_1^0,x_2^0,\cdots,x_m^0)^T,
  \lbl{4.1}\ees
then the particle system is called $N$-periodic system in
distributional sense.
  \end{defi}
One can give some examples to satisfy (\ref{4.1}). For example, suppose a particle system
has five states and the initial distribution is
$(\frac{1}{10},\frac{1}{10},\frac{7}{20},\frac{2}{5},\frac{1}{20})^T$. Assume that the one-step
transition probability matrix is
 \bess
P=\left(\begin{array}{cccccccc}
 a_1 \ &\ \ b_1 &\ \ 0 &\ \  0  &\ \ 0  \\
 a_2 &\ \ b_2 &\ \ 0 &\ \  0  &\ \ 0  \\
 0 \ &\ \ 0 &\ \ 0 &\ \  1 &\ \ 0 \\
 0 \ &\ \ 0 &\ \ 0&\ \ 0 &\ \ 1 \\
 0 \ &\ \ 0 &\ \ 1 &\ \  0 &\ \ 0
\end{array}\right), \nonumber
   \eess
where $0\leq a_i,b_i\leq1$, $a_i+b_i=1$, $i=1,2$. Then one can find that
   \bess
P^2\left(\frac{1}{10},\frac{1}{10},\frac{7}{20},\frac{2}{5},\frac{1}{20}\right)^T
=\left(\frac{1}{10},\frac{1}{10},\frac{7}{20},\frac{2}{5},\frac{1}{20}\right)^T.
   \eess
On the other hand, it is easy to see that if
   \bes
P^{(N)}=I_m,
   \lbl{4.2}\ees
where $I_m$ denotes $m\times m$ identity matrix (which is called Idempotent matrix in
algebra), then the equality (\ref{4.1}) holds.
A stochastic process is called \emph{strong $N$-periodic system in
distributional sense} if (\ref{4.2}) holds.
We remark that the number $N$ in (\ref{4.2}) is definitely equal to
least common multiple of the periodicity of every particle.

For continuous time and continuous state stochastic process, we borrow the idea of
\cite{BKRS,Chenfeng2017}. A stochastic process is called $T$-periodic system in
distributional sense if $\mu(t+T,x)=\mu(t,x)$ for all $t\geq0$ and $x\in\mathbb{R}^d$,
where $\mu$ is defined as in Proposition \ref{p1.1}.

Before we close this section, we establish the existence of periodic
solution in distributional sense on bounded domain. We generalize the result of
\cite{Chenfeng2017} to the bounded domain. We remark that
the boundary of the bounded domain should be reflective. If the boundary is absorbing, then we cannot get the
limit in the following sense
   \bess
\mu_n(f)=\int_Dfd\mu_n\to\int_Dfd\mu=\mu(f),\ \ \ as\ \ n\to\infty,
  \eess
where $\mu_n$ and $\mu$ are probability measure of some stochastic
process on the bounded domain $D\subset\mathbb{R}^d$. The probability
measure considered here keeps entirety, i.e., $\mu_n(\bar D)=1$ and the limit probability measure $\mu_0(\bar D)=1$.
The results obtained here coincide with those in next section.

Let $D$ be a convex domain in $\mathbb{R}^d$ and $(\Omega,\mathcal{F},P)$ be a complete probability
space with an increasing family $\{\mathcal{F}_t\}_{t\geq0}$
of sub-$\sigma$-fields of $\mathcal{F}$. Suppose an $\mathcal{F}_t$-adapt
$r$-dimensional Brownian motion $B(t)=(B^1(t),\cdots,B^r(t))$
with $B(0)=0$ is given. Let $\sigma(t,x)$ and $b(t,x)$ be
$\mathbb{R}^d\otimes\mathbb{R}^d$-valued and $\mathbb{R}^d$-valued
functions, both being defined on $\mathbb{R}_+\times\bar D$, respectively. Consider the stochastic differential equation with
reflection
   \bes
dX_t=b(t,X)dt+\sigma(t,X)dB+d\Phi,\ \ \ X(0)=x,
  \lbl{4.3}\ees
where $x\in\bar D$ and $\{\Phi(t)\}$ is an associated process of
$\{X(t)\}$. In \cite{LionsSznitman1984}, the authors gave the
relationship between $\Phi$ and $X$, i.e.,
   \bess
\Phi_t=\int_0^t\nu(X_s)d|\Phi|_s,\ \ \ |\Phi|_t=\int_0^t1_{\{X_s\in\partial D\}}d|\Phi|_s,
   \eess
where $\nu$ is the unit outward normal to $\partial D$ at $x$, and
$k_t$ stands for the total variation of $k$ on $[0,t]$.
In order to make the meaning of $\Phi_t$ clearly, we introduce
the following spaces of functions, see \cite[Page 164]{Tanaka1979} for more details.

$C(\mathbb{R}_+,\mathbb{R}^d)$ (resp. $C(\mathbb{R}_+,\bar D)$)
$=$ the space of $\mathbb{R}^d$-valued (resp. $\bar D$-valued)
continuous functions on $\mathbb{R}_+$.

$\mathbb{D}(\mathbb{R}_+,\mathbb{R}^d)$ (resp. $\mathbb{D}(\mathbb{R}_+,\bar D)$)
$=$ the space of $\mathbb{R}^d$-valued (resp. $\bar D$-valued)
right continuous functions on $\mathbb{R}_+$ with left limits.

On $C(\mathbb{R}_+,\mathbb{R}^d)$ and $C(\mathbb{R}_+,\bar D)$
we consider the compact uniform topology. Given a function $\xi$
in $\mathbb{D}(\mathbb{R}_+,\bar D)$, a function $\Phi$ is said to be
associated with $\xi$ if the following three conditions are
satisfied.

(i) $\Phi$ is a function in $\mathbb{D}(\mathbb{R}_+,\mathbb{R}^d)$ with
bounded variation and $\Phi(0)=0$.

(ii) The set $\{t\in\mathbb{R}_+:\xi(t)\in D\}$ has $d|\Phi|$-measure zero.

(iii) For any $\eta\in C(\mathbb{R}_+,\bar D)$, $(\eta(t)-\xi(t),\Phi(dt))\geq0$.

Using the above properties, Tanaka proved that the following
Lemma.

\begin{prop}\lbl{p4.2}\cite[Lemma 2.2]{Tanaka1979} Let $w,\tilde w\in \mathbb{D}(\mathbb{R}_+,\mathbb{R}^d)$ with $w(0),\tilde w(0)\in \bar D$,
and $\xi,\tilde\xi$ be any solutions of
   \bess
\xi=w+\Phi,\ \ \ \ \tilde\xi=\tilde w+\tilde\Phi,
 \eess
respectively. Then we have
   \bess
|\xi(t)-\tilde\xi(t)|^2&\leq&|w(t)-\tilde w(t)|^2\\
&&+2\int_0^t(w(t)-\tilde w(t)-w(s)+\tilde w(s),\Phi(ds)-\tilde\Phi(ds)).
   \eess
   \end{prop}

Tanaka \cite{Tanaka1979} obtained the following result.
\begin{prop}\cite[Theorem 4.1]{Tanaka1979}\lbl{p4.1} If there exists a constant
$K>0$ such that
   \bess
&&\|\sigma(t,x)-\sigma(t,y)\|\leq K|x-y|,\ \ \|b(t,x)-b(t,y)\|\leq K|x-y|, \\
&&\|\sigma(t,x)\|\leq K(1+|x|^2)^{1/2},\ \ \ \|b(t,x)\|\leq K(1+|x|^2)^{1/2},
  \eess
then there exists a (pathwise) unique $\mathcal{F}_t$-adapted solution
(\ref{4.3}) for any $x\in\bar D$.
\end{prop}
Later, Lions-Sznitman \cite{LionsSznitman1984} generalized the
results of \cite{Tanaka1979}. Now, we follow the idea of
\cite{Chenfeng2017} to prove the existence of periodic solution
in distributional sense on bounded domain.
Note that the bounded domain with reflection boundary is similar to
the whole space, the proof is similar to that of \cite{Chenfeng2017}.
We only write out the difference. Due to that nothing is lost in the
bounded domain, so the probability measure  on $\bar D$
will be always $1$. Using this fact, we can obtain a similar theorem
on bounded domain to \cite[Theorem Page 9]{Skorokhod1965}. And thus
Lemmas 2.3 and 2.4 in \cite{Chenfeng2017} hold for the bounded domain.

Let $\mathcal{P}(\bar D)$
be the set of Borel probability measures on $\bar D$. We denote
the law of $X$ on $\bar D$ by $\mu:\,\mathbb{R}\to\mathcal{P}(\bar D)$.
Assume there exists a stochastic process $L$ such that the solution
$Y(t)$ on $\mathbb{R}_+$ of (\ref{4.3}) satisfying
   \bes
|Y(nT)|\leq L, \ \ \ \mathbb{E}|L|^2<\infty, \ \ \ n=1,2,\cdots,.
  \lbl{4.4}\ees
We borrow symbols from the \cite{Chenfeng2017}. $P\circ[Y(t)]^{-1}$
denotes the distribution of $Y(t)$.  Similar to Section 2 of \cite{Chenfeng2017}, we define the $d_{BL}$ which
means the distance of
bounded and Lipschitz function.
   \bess
&&\|h\|_\infty=\sup_{\bar D}|h(x)|,\ \
\|h\|_L=\sup_{x,y\in \bar D,x\neq y}\{\frac{|h(x)-h(y)|}{|x-y|}\},\\
&&\|h\|_{BL}=\max\{\|h\|_\infty,\,\|h\|_L\},\ \
d_{BL}(\mu,\nu)=\sup_{\|h\|_{BL}\leq1}\Big|\int hd(\mu-\nu)\Big|
   \eess
for all $\mu,\nu\in \mathcal{P}(\bar D)$ and all Lipschitz
continuous real-valued functions $h$ on $\bar D$.
It is easy to check that $(d_{BL},\mathcal{P}(\bar D))$ is
a complete metric space, see \cite[Page 390]{Dudley2002} for details.
The main result is the following theorem.
\begin{theo}\lbl{t4.1} Let $b$ and $\sigma$ be continuous functions
which are $T$-periodic in the time variable and satisfy the assumptions
of Proposition \ref{p4.1}. If (\ref{4.4}) holds  and
  \bes
\lim\limits_{k\to\infty}
\frac{1}{n_k+1}\sum_{m=0}^{n_k}
d_{BL}(P\circ[Y((m+1)T)1_{A_m}]^{-1},P\circ[Y(mT)1_{A_m}]^{-1})=0,
   \lbl{4.5}\ees
where $Y(t)$ is a solution of (\ref{4.3}), $A_m$ is defined as in (\ref{4.01}) and $\{n_k\}$ is a sequence of
integers tending to $+\infty$ and $d_{BL}$ is a metric,
then there
exists an $L^2$-bounded $T$-periodic solution in distribution sense of
(\ref{4.3}).
  \end{theo}

{\bf Proof.} Inspired by \cite{Chenfeng2017}, define the stochastic
process
   \bess
X_k(0,\omega)=Y(\chi_k(\omega)), \ \ \
X_k(t,\omega)=Y(t+\chi_k(\omega),\omega),
   \eess
where $\omega\in\Omega$, $\chi_k$ is a random variable independent of
$B_t$ and $Y(0,\omega)$ such that $P(\chi_k=nT)=\frac{1}{k+1}$, $n=0,1,\cdots,k$. Due to the functions $b$ and $\sigma$ are
$T$-periodic in time variable and the fact that
$\Phi_t$ just depends on $X_t$, $X_k$ is still a solution of
   \bess
dX_t=b(t,X_t)dt+\sigma(t,X_t)d\tilde B_t+d\Phi_t,
\ \ \ \tilde B_t=B(t+nT)-B(nT),
   \eess
where $\tilde B_t$ has the same distribution with $B_t$.
Similar to \cite{Chenfeng2017}, using the fact that $\chi_k$ is
independent of $\tilde B_t$, we have
   \bes
P(X_k(t)\in A_0)=\frac{1}{k+1}\sum_{n=0}^kP(Y(t+nT)\in A_0),
  \lbl{4.6} \ees
where $A_0\subset\bar D$ is a Borel set. It follows from
(\ref{4.4}), (\ref{4.6}) and Chebyshev's inequality that
  \bess
P(|X_k(0,\omega)|>R)\leq\frac{1}{k+1}\sum_{n=0}^k
\frac{\mathbb{E}|Y(nT)|^2}{R^2}\leq\frac{C}{R^2}\to0\ \ {\rm as}\ R\to\infty.
  \eess
Applying Skorokhod' Lemma (\cite[Lemmas 2.3 and 2.4]{Chenfeng2017}),
we have that in some probability space
$(\tilde\Omega,\mathcal{\tilde F},\tilde P)$ there exists a sequence $\tilde X_k(0,\tilde\omega)$ ($k=0,1,\cdots$) with the same
distribution as $X_k(0,\omega)$ such that some subsequence
$\{\tilde X_{n_k}(0,\tilde\omega)\}_{k=0,1,\cdots}$ converges in probability
to $\tilde X(0,\tilde\omega)$. Also, we can construct random variables
$X_k(\omega)$ and $X(\omega)$ on the space $(\Omega,\mathcal{F},P)$, whose joint distribution is the same
as the joint distribution of $\tilde X_k(\tilde\omega)$ and
$\tilde X(\tilde\omega)$. Notice that $\tilde X_{n_k}(0,\tilde\omega)$ has the same distribution as $ X_{n_k}(0,\omega)$,
and thus we have $\mathbb{E}|\tilde X_{n_k}(0,\tilde\omega)|^2=\mathbb{E}|X_{n_k}(0,\omega)|^2$,
$|\tilde X_{n_k}(0,\tilde\omega)|\leq L$ and $\mathbb{E}|L|^2<\infty$. The Vitali's theorem implies that
$\mathbb{E}|\tilde X_{n_k}(0,\tilde\omega)-\tilde X(0,\tilde \omega)|^2\to0$ as $k\to\infty$.

Let $\tilde X_{n_k}(t)$ be the solution of
   \bess
dX_t=b(t,X_t)dt+\sigma(t,X_t)d\tilde B_t+d\Phi_t
  \eess
with initial data $\tilde X_{n_k}(0,\tilde\omega)=X_{n_k}(\tilde \omega)$ on the probability space
$(\tilde\Omega,\mathcal{\tilde F},\tilde P)$.
Note that
   \bess
\tilde X_{n_k}(t)-\tilde X(t)&=&\tilde X_{n_k}(\tilde\omega)-\tilde X(\tilde\omega)+\int_0^t(b(s,\tilde X_{n_k}(t))-b(s,\tilde X(s)))ds\\
&&+\int_0^t(\sigma(s,\tilde X_{n_k}(t))-\sigma(s,\tilde X(s)))d\tilde B_s+\Phi_t-\tilde\Phi_t,
  \eess
we have
   \bess
|\tilde X_{n_k}(t)-\tilde X(t)|^2&\leq&3|\tilde X_{n_k}(\tilde\omega)-\tilde X(\tilde\omega)|^2+3|\int_0^t(b(s,\tilde X_{n_k}(t))-b(s,\tilde X(s)))ds|^2\\
&&+3|\int_0^t(\sigma(s,\tilde X_{n_k}(t))-\sigma(s,\tilde X(s)))d\tilde B_s+\Phi_t-\tilde\Phi_t|^2.
  \eess
Let
   \bess
&&\xi(t)=w(t)+\Phi_t,\ \ \tilde\xi(t)=\tilde w(t)+\tilde\Phi_t,\\
&&w(t)=\int_0^t\sigma(s,\tilde X_{n_k}(t))d\tilde B_s, \ \ \
\tilde w(t)=\int_0^t\sigma(s,\tilde X(t))d\tilde B_s.
  \eess
Then applying Proposition \ref{p4.2}, using It\^{o}
isometry, and noting that $w(t)$ is a
martingale with respect to $\mathcal{\tilde F}$, we have
(also see the proof of \cite[Theorem 4.1]{Tanaka1979}, and here
we use the reason why "the reminder" disappear in (4.4) on Page
175)
   \bess
\mathbb{E}|\tilde X_{n_k}(t)-\tilde X(t)|^2&\leq&3\mathbb{E}|\tilde X_{n_k}(\tilde\omega)-\tilde X(\tilde\omega)|^2+3t\mathbb{E}\int_0^t\|b(s,\tilde X_{n_k}(s))-b(s,\tilde X(s))\|^2ds\\
&&+3\mathbb{E}\big|\int_0^t(\sigma(s,\tilde X_{n_k}(t))-\sigma(s,\tilde X(s)))d\tilde B_s\big|^2\\
&&+6\mathbb{E}\int_0^t(w(t)-\tilde w(t)-w(s)+\tilde w(s),\Phi(ds)-\tilde\Phi(ds))\\
&\leq&3\mathbb{E}|\tilde X_{n_k}(\tilde\omega)-\tilde X(\tilde\omega)|^2+3tK^2\mathbb{E}\int_0^t|\tilde X_{n_k}(s))-\tilde X(s)|^2ds\\
&&+6\mathbb{E}\int_0^t|\sigma(s,\tilde X_{n_k}(t))-\sigma(s,\tilde X(s))|^2ds\\
&\leq&3\mathbb{E}|\tilde X_{n_k}(\tilde\omega)-\tilde X(\tilde\omega)|^2+3K^2(t+2)\mathbb{E}\int_0^t|\tilde X_{n_k}(s))-\tilde X(s)|^2ds,
  \eess
where we used the fact that (independent increment of Brownian motion)
   \bess
\mathbb{E}\int_0^t(w(t)-\tilde w(t)-w(s)+\tilde w(s),\Phi(ds)-\tilde\Phi(ds))=0.
  \eess
By Gronwall's inequality, we have
   \bess
\mathbb{E}|\tilde X_{n_k}(t)-\tilde X(t)|^2
\leq 3\mathbb{E}|\tilde X_{n_k}(\tilde\omega)-\tilde X(\tilde\omega)|^2e^{3(t+2)K^2}\to0\ \
{\rm as }\ k\to\infty.
  \eess
Since the uniqueness of weak solutions implies the
uniqueness of laws, we have
   \bess
P\circ[X_{n_k}(t)]^{-1}=P\circ[\tilde X_{n_k}(t)]^{-1}
\to P\circ[\tilde X(t)]^{-1}
   \eess
unformly on $[0,T]$. Moreover, we can replace $\tilde X(0,\tilde \omega)$ on $(\tilde\Omega,\mathcal{\tilde F},\tilde P)$
by $X(0,\tilde \omega)$ on $(\Omega,\mathcal{F},P)$ with the
same law. Then the solution $X(t)$ admits the same distribution
of $\tilde X(t)$ by weak uniqueness of the equation (\ref{4.3}).
It suffices to prove that
   \bess
P\circ[X(T)]^{-1}=P\circ[X(0)]^{-1}.
  \eess
Denote
   \bes
A_m=\{\omega\in\Omega:\ \  \chi_{n_k}(\omega)=mT\},\ \ m=0,1,\cdots,n_k,
   \lbl{4.01}\ees
then we have
   \bess
P(A_m)=\frac{1}{n_k+1}, \ \ m=0,1,\cdots,n_k.
  \eess
 By using the above equality, we get
    \bes
\int_\Omega\phi(X_{n_k}(T)) dP=\frac{1}{n_k+1}\sum_{m=0}^{n_k}\int_{A_m}\phi(Y((m+1)T)) dP.
   \lbl{4.7}\ees
Set
   \bess
1_{A_m}(\omega)=\left\{\begin{array}{lll}
1,\ \ \ \omega\in A_m,\\
0,\ \ \ \omega\not\in A_m.
 \end{array}\right.\eess
It follows from (\ref{4.5}), (\ref{4.6}) and (\ref{4.7}) that (see the proof of \cite{Chenfeng2017} in Page 292 for details)
   \bes
&&d_{BL}(P\circ[X(T)]^{-1},P\circ[X(0)]^{-1})\nonumber\\
&=&\lim\limits_{k\to\infty}d_{BL}
(P\circ[X_{n_k}(T)]^{-1},P\circ[X_{n_k}(0)]^{-1})\ \ \ \ \ {\rm (by\ the\ definition)}\nonumber\\
&=&\lim\limits_{k\to\infty}\sup_{\|\phi\|_{BL}\leq1}
\left|\int_D\phi dP\circ[X_{n_k}(T)]^{-1}-\int_D\phi d P\circ[X_{n_k}(0)]^{-1}\right|\nonumber\\
&=&\lim\limits_{k\to\infty}\sup_{\|\phi\|_{BL}\leq1}
\left|\int_\Omega\phi(X_{n_k}(T)) dP-\int_\Omega\phi(X_{n_k}(0)) d P\right|\nonumber\\
&=&\lim\limits_{k\to\infty}\sup_{\|\phi\|_{BL}\leq1}
\left|\frac{1}{n_k+1}\sum_{m=0}^{n_k}\int_{A_m}\phi(Y((m+1)T)) dP-\frac{1}{n_k+1}\sum_{m=0}^{n_k}\int_{A_m}\phi(Y(mT)) d P\right|\nonumber\\
&=&\lim\limits_{k\to\infty}\sup_{\|\phi\|_{BL}\leq1}\left|
\frac{1}{n_k+1}\sum_{m=0}^{n_k}\int_{A_m}[\phi(Y((m+1)T)) -\phi(Y(mT)) ]d P\right|\nonumber\\
&\leq&\lim\limits_{k\to\infty}\sup_{\|\phi\|_{BL}\leq1}
\frac{1}{n_k+1}\sum_{m=0}^{n_k}\left|\int_{A_m}[\phi(Y((m+1)T)) -\phi(Y(mT)) ]d P\right|\nonumber\\
&\leq&\lim\limits_{k\to\infty}
\frac{1}{n_k+1}\sum_{m=0}^{n_k}
d_{BL}(P\circ[Y((m+1)T)1_{A_m}]^{-1},P\circ[Y(mT)1_{A_m}]^{-1})\nonumber\\
&=&0,
  \lbl{4.8}\ees
that is to say, $X(T)$ has the same distribution as $X(0)$.

Define the function $z:\mathbb{R}_+\to\bar D$ by
   \bess
z(t)=Y(t-n_tT),
  \eess
where $n_t=\max\{n\in \mathbb{N}|nT\leq t\}$. Then $z(t)$ is a
$T$-periodic solution to (\ref{4.3}). The proof is complete.
$\Box$

\begin{remark}\lbl{r4.1}
Comparing with the assumptions in Theorem \ref{t4.1} with \cite[Theorem 1.2]{Chenfeng2017},
one can find there is a little difference from \cite[Theorem 1.2]{Chenfeng2017}.
The reason is that in (\ref{4.8}), the second last equality we use an
equality, which different from \cite{Chenfeng2017}, where they used the following inequality
  \bess
&&\lim\limits_{k\to\infty}\sup_{\|\phi\|_{BL}\leq1}
\frac{1}{n_k+1}\sum_{m=0}^{n_k}\left|\int_{A_m}[\phi(Y((m+1)T)) -\phi(Y(mT)) ]d P\right|\nonumber\\
&\leq&\lim\limits_{k\to\infty}\sup_{\|\phi\|_{BL}\leq1}
\frac{1}{n_k+1}\sum_{m=0}^{n_k}d_{BL}(P\circ[Y((m+1)T)]^{-1},P\circ[Y(mT)]^{-1}).
  \eess
Noting that
  \bess
d_{BL}(P\circ[Y((m+1)T)1_{A_m}]^{-1},P\circ[Y(mT)1_{A_m}]^{-1})\leq d_{BL}(P\circ[Y((m+1)T)]^{-1},P\circ[Y(mT)]^{-1}),
   \eess
implies that the condition (\ref{4.5}) is weaker than \cite[(5)]{Chenfeng2017}.

On the other hand, if $\sigma=0$, the condition (\ref{4.5}) becomes
   \bess
\lim\limits_{k\to\infty}\sup_{\|\phi\|_{BL}\leq1}
\frac{1}{n_k+1}\sum_{m=0}^{n_k}1_{A_m}
|Y((n_k+1)T)-Y(0)|=0.
   \eess
which is an extension assumption to Halanay \cite{Halanay1966}.

Similarly, if we define
   \bess
P(A_m)=p_m,
   \eess
then the condition (\ref{4.5}) becomes
  \bess
\lim\limits_{k\to\infty}
\sum_{m=0}^{n_k}p_m
d_{BL}(P\circ[Y((m+1)T)1_{A_m}]^{-1},P\circ[Y(mT)1_{A_m}]^{-1})=0.
   \eess
For different choice of $p_m$, we can get the different periodic
solution, and thus there are infinity periodic solutions in
distributional sense.
\end{remark}

One can also obtain the existence of periodic solution to
(\ref{3.11}) similar to \cite[Theorem 1.1]{Chenfeng2017}.
Base on the Theorem \ref{t4.1}, one can establish the
uniqueness of periodic solution under similar assumptions
to \cite[Theorem 1.4]{Chenfeng2017} and give the Lyapunov
functional to verify the existence of periodic solution
to (\ref{4.3}), see \cite[Theorem1.3]{Chenfeng2017}.
Due to similarity, we omit the details to readers.

\section{Existence of periodic solutions on bounded domain and in the whole space}\label{sec4}\setcounter{equation}{0}
In this section, we obtain some properties of PDF by
considering the existence of periodic solutions to the Fokker-Planck equations.
In order to establish the desired results, we divide this section into two parts.
\subsection{Bounded domain with Dirichlet boundary condition}
The reason why we first consider the Dirichlet boundary condition problem is the
assumptions on $(\mathcal{A},\mathcal{B})$. More precisely, in Section 2, we make
the assumption that $\mathcal{B}=\mathcal{B}(x,D)$ is independent of $t\in[0,T]$.
In this subsection, we assume that a stochastic process $X_t$ satisfies
   \bes\left\{\begin{array}{llll}
dX_t=b(t,X_t)dt+\sigma(t,X_t)dB_t,\\
 X_0=x,
   \end{array}\right.\lbl{3.1}\ees
where $b$ is an $d$-dimensional vector function, $\sigma$ is an $d\times m$ matrix function and
$B_t$ is an $m$-dimensional Brownian motion.
Our aim is to study the properties of PDF by considering the  existence of periodic solution to the
corresponding Fokker-Planck equation. Throughout this subsection, we assume the
particle (in the system) will die if it touch the boundary. That is to say, the PDF of this system
satisfies the following evolution equation
   \bes\left\{\begin{array}{llll}
  \partial_tp(t,x)=\mathcal{A}^*p(t,x), \ \ &in\  D\times(0,T],\\
p=0, \ \ \ \ &on\ \partial D\times(0,T],\\
p(0,x)=p_0(x), \ \ \ & in \ D,
   \end{array}\right.\lbl{3.2}\ees
where
   \bess
\mathcal{A}^*p=-\sum_i\frac{\partial}{\partial x_i}(b_ip)+\frac{1}{2}\sum_{i,j}
\frac{\partial^2}{\partial x_i\partial x_j}((\sigma\sigma^T)_{ij}p).
  \eess
In order to get the properties of $p$ in (\ref{3.2}), we need consider the following auxiliary equation
   \bes\left\{\begin{array}{llll}
    \partial_tu(t,x)=\mathcal{A}^*u(t,x)+\mu u, \ \ &in\  D\times(0,T],\\
u=0, \ \ \ \ &on\ \partial D\times(0,T],\\
u(0,x)=u(T,x), \ \ \ & in \ D.
   \end{array}\right.\lbl{3.3}\ees
 In order to get the existence of solution of (\ref{3.3}), we first consider the following initial boundary problem
   \bes\left\{\begin{array}{llll}
  \partial_tu(t,x)=\mathcal{A}^*u(t,x)+\mu u(t,x), \ \ &in\  D\times(0,T],\\
u=0, \ \ \ \ &on\ \partial D\times(0,T],\\
u(0,x)=u_0(x), \ \ \ & in \ D,
   \end{array}\right.\lbl{3.4}\ees
where $u_0(x)$ is a fixed function and will be given later.
It is easy to see that there exists an evolution operator $U(t,s)$ such that the solution of (\ref{3.4}) can be represented
in the following form (see section 13 in \cite{Hess1991book})
   \bess
u(t)=U(t,0)u_0+\mu\int_0^tU(t,s)u(s)ds.
   \eess
If we set $v(t,x)=e^{-\mu t}u(t,x)$, then we have
   \bess
 v(t,x)=U(t,0)v_0=U(t,0)u_0,
   \eess
that is to say,
   \bess
u(t,x)=e^{\mu t}U(t,0)u_0.
   \eess
It is easy to check that $v(t,x)$ is a solution of (\ref{3.2}) with $p_0=u_0$.

We assume that

(C2)  The operator $\mathcal{A}^*$ is a uniform elliptic operator, $(\sigma\sigma^T)_{ij}(t,x)\in C^{2+\theta,1+\frac{\theta}{2}}(\bar D\times\mathbb{R})$,
$b_i(t,x)\in C^{1+\theta,\frac{1+\theta}{2}}(\bar D\times\mathbb{R})$, and both $(\sigma\sigma^T)_{ij}(t,x)$
and $b_i(t,x)$ are T-periodic  in $ t$.

Set $K=U(T,0)$. Then the Krein-Rutman theorem \cite[Theorem 7.2]{Hess1991book} implies that $r:=spr(K)>0$.
Proposition \ref{p2.3} yields that $\mu>0$ if $r<1$.

\begin{lem}\lbl{l3.1} Assume that the condition (C2) holds. Assume further that
  \bess
a_0(t,x):=-\frac{1}{2}\sum_{i,j}\frac{\partial^2}{\partial x_i\partial x_j}(\sigma\sigma^T)_{ij}+div b(t,x)\geq0.
  \eess
Then $0<r<1$, where $r:=spr(K)$.
\end{lem}

{\bf Proof.} For completely, we give the outline of the proof. To show that $r<1$, let $u_0\in W^{2,p}_{0}(D)$,
$u_0\gg0$ be a principal eigenfunction of $K$, i.e., $Ku_0=ru_0$.
Then $u:=U(\cdot,0)u_0$ solves
   \bess\left\{\begin{array}{llll}
  \partial_tu(t,x)-\mathcal{A}^*u(t,x)=0, \ \ & in\  D\times(0,T],\\
u=0, \ \ \ \ &on\ \partial D\times(0,T],\\
u(0,x)=u_0(x), \ \ \ & in \ D.
   \end{array}\right.\eess
If $a_0\geq0$ in $\bar D\times[0,T]$ and $u_0>0$ in $W^{2,p}_{0}(D)$, the Propositions
13.1 and 13.3 and Remark 13.2 in \cite{Hess1991book} that $U(t,\tau)u_0\gg0$ in $W^{2,p}_{0}(D)$
for $\tau<t\leq T$. Now let $v:=\|u_0\|_{C(\bar D)}-u$, then $v$ satisfies
   \bes\left\{\begin{array}{llll}
  \partial_tv(t,x)-\mathcal{A}^*v(t,x)=a_0\|u_0\|_{C(\bar D)}\geq0, \ \ & in\  D\times(0,T],\\
v\geq0, \ \ \ \ &on\ \partial D\times(0,T],\\
v(0,x)\geq0, \ \ \ & in \ D,
   \end{array}\right.\lbl{3.5}\ees
and hence $v\gg0$ in $W^{2,p}_{0}(D)$ for each $0<t\leq T$ by
the Propositions 13.1 and 13.3 and Remark 13.2 in \cite{Hess1991book}. In particular,
   \bess
 r\|u_0\|_{C(\bar D)}=\|Ku_0\|_{C(\bar D)}=\|u(T)\|_{C(\bar D)}<\|u_0\|_{C(\bar D)},
   \eess
which implies that $r<1$. $\Box$

\begin{theo}\lbl{t3.1} Under the assumptions of Lemma \ref{l3.1}, there exists a unique solution to the
equation (\ref{3.3}). Hence the solution of (\ref{3.2}) satisfies exponential decay for any fixed point
in $D$ as time goes to infinity under the special initial data. That is to say, the solution $p(t,x)$
has the following property:
   \bess
p(nT,x)=e^{-\mu nT}p_0(x), \ \ x\in D,\ n\in Z,
   \eess
where $\mu>0$ is given as in (\ref{3.3}) and $p_0(x)$ satisfies $Kp_0=rp_0$ with $K=U(T,0)$.
\end{theo}

{\bf Proof.} It follows from Lemma \ref{l3.1} that the principal eigenvalue of $K$ satisfies
$0<r<1$. Proposition \ref{p2.3} implies that
    \bess
\mu=-\frac{1}{T}\log r
   \eess
is an eigenvalue of (\ref{3.3}) with the positive eigenfunction
   \bess
u(t)=e^{\mu t}U(t,0)u_0.
   \eess
Take $u_0$ be the principal eigenfunction of $K$ (take $p_0(x)$ be the principal eigenfunction of corresponding $K$ in
equation (\ref{3.3})). And thus we have
   \bess
u(T)=e^{\mu T}U(T,0)u_0=\frac{1}{r}Ku_0=u(0).
  \eess
That is to say, $u(t)$ is the solution of (\ref{3.3}). The uniqueness of principal eigenvalue implies
the uniqueness of solution of (\ref{3.3}). By $T$-periodicity of $\mathcal{A}(t)$, we have $U(t,\tau)=U(t+nT,\tau+nT)$,
$n\in Z$. Noting that the solution of (\ref{3.2}) can be written as
   \bess
p(t,x)=e^{-\mu t}u(t,x),
  \eess
we have
   \bess
p(T,x)=e^{-\mu T}u(T,x)=e^{-\mu T}u(0,x)=e^{-\mu T}p_0(x).
  \eess
By using the properties
\bess
U(t,t)=I, \ \ \ U(s,t)U(t,\tau)=U(s,\tau), \ \ 0\leq\tau\leq t\leq s\leq T,
  \eess
we have
   \bess
p(nT,x)=U(nT,0)p_0(x)=U(nT,(n-1)T)\cdots U(T,0)p_0(x)=e^{-\mu nT}p_0(x),
   \eess
where we used the fact that $U(T,0)p_0(x)=e^{-\mu T}p_0(x)$. $\Box$

\begin{remark}\lbl{r3.1} In the proof of Lemma \ref{l3.1}, we know that
the initial data $u_0$ (or $P_0$) is special function, that is, $u_0$ satisfies
$Ku_0=ru_0$. Now, we give an example to show that this is possible.
Consider the problem
   \bess\left\{\begin{array}{llll}
  \partial_tu(t,x)-\Delta u(t,x)=0, \ \ & in\  D\times(0,T],\\
u=0, \ \ \ \ &on\ \partial D\times(0,T],\\
u(0,x)=u_0(x), \ \ \ & in \ D.
   \end{array}\right.\eess
Assume that $u_0$ satisfies ($r>0$)
     \bess
-\Delta u_0=ru_0, \ \ in \ D,\ \ \ \ u|_{\partial D}=0,
   \eess
and by using the fact
   \bess
u(t,x)=e^{t\Delta_D}u_0,
  \eess
where $\Delta_D$ denotes the Laplace operator with Dirichlet boundary, then we get (by Taylor expansion)
   \bess
u(t,x)=\sum_{i=0}^\infty \frac{(t\Delta_D)^i}{i!}u_0=\sum_{i=0}^\infty \frac{(-tr)^i}{i!}u_0=e^{-tr}u_0.
  \eess
And thus the solution of the following equation
   \bess\left\{\begin{array}{llll}
  \partial_tv(t,x)-\Delta v(t,x)=\mu v(t,x), \ \ & in\  D\times(0,T],\\
v=0, \ \ \ \ &on\ \partial D\times(0,T],\\
v(0,x)=u_0(x), \ \ \ & in \ D.
   \end{array}\right.\eess
can be written as
   \bess
v(t,x)=e^{\mu t}u(t,x)=e^{t(\mu-r)}u_0.
   \eess
If we want to get $v(T,x)=v_0(x)$, we will take $\mu=r$. Because
there is no concrete value for $T$, we obtain that the solution $u$
satisfies $u(t,x)=e^{-\mu t}u_0$ for all $t>0$ and $x\in D$.
  \end{remark}
It follows from Theorem \ref{t3.1} that it is difficult to obtain the existence of periodic solution to linear
parabolic equation. Now we turn to the nonlinear case.
In 1998 Pardoux and Zhang proved in \cite{PardouxZhang1998} a probabilistic
formula for the viscosity solution of a system of semilinear PDEs with Neumann
boundary condition
   \bess
 \left\{\begin{array}{llll}
  \partial_tu+\mathcal{A}(u)+f(t,x,u)(t,x))=0, \ \ & in\  D\times(0,T],\\
\frac{\partial u}{\partial \nu}+g(t,x,u)=0, \ \ \ \ &on\ \partial D\times(0,T],\\
u(T,x)=h(x), \ \ \ & in \ D,
   \end{array}\right.\eess
where $D$ is an open connected bounded subset of $\mathbb{R}^d$.

In order to get the existence of periodic solution for Dirichlet
problem, we need to consider the nonlinear parabolic equation
   \bes\left\{\begin{array}{llll}
  \partial_tu(t,x)=\mathcal{A}^*u(t,x)+f(t,x,u), \ \ &in\  D\times(0,T],\\
u=0, \ \ \ \ &on\ \partial D\times(0,T],\\
u(0,x)=u_0(x), \ \ \ & in \ D.
   \end{array}\right.\lbl{3.6}\ees
The assumptions on $f$ will be given later.

We first recall some results.
In \cite{Hess1991book}, Hess considered the following periodic initial boundary problem
       \bes\left\{\begin{array}{llll}
  \partial_tu(t,x)+\mathcal{A}(t)u(t,x)=f(x,t,u), \ \ &in\  D\times(0,T],\\
\mathscr{B}u:=\partial_\nu u+bu=0, \ \ \ \ &on\ \partial D\times(0,T],\\
u(0,x)=u(T,x), \ \ \ & in \ D,
   \end{array}\right.\lbl{3.8}\ees
where they assumed the function $b$ does not depend on $t$, and
   \bess
\mathcal{A}(t)u=-\sum_{i,j}a_{ij}(t,x)\frac{\partial^2}{\partial{x_j}\partial{x_k}}u+\sum_ia_i(t,x)
\frac{\partial}{\partial{x_i}}u+a_0(t,x)u.
   \eess
They used the upper and lower solution method to prove the existence of periodic solution of
(\ref{3.8}). We first recall the definition of upper (lower) solution.
  \begin{defi}\lbl{d3.1} Let $u\in C^{1,0}([0,T]\times\bar D)\cap C^{2,1}([0,T)\times D)$.
  Such a function $u$ is referred to as an upper (lower) solution if
       \bess\left\{\begin{array}{llll}
  \partial_tu(t,x)+\mathcal{A}(t)u(t,x)\geq(\leq)f(x,t,u), \ \ &in\  D\times(0,T],\\
\mathscr{B}u\geq(\leq)0, \ \ \ \ &on\ \partial D\times(0,T],\\
u(0,x)\geq(\leq)u(T,x), \ \ \ & in \ D,
   \end{array}\right.\eess
\end{defi}
Let $\bar u\geq\underline u$ be the upper and lower solution of (\ref{3.8}), respectively.
We define $\sigma=\min_{[0,T]\times\bar D}\underline u$ and $\omega=\max_{[0,T]\times\bar D}\bar u$.
Set
   \bess
W_{p,\mathscr{B}}^{2}(D)=\{u\in W^2_p(D): \ [\partial_\nu u+b(x)u=0\}, \ \ \ p>d.
   \eess
\begin{prop}\lbl{p3.1} Suppose $\bar u\geq\underline u$ are the upper and lower solutions of (\ref{3.8}), respectively.
Let $f(\cdot,\cdot,u)\in C^{\alpha/2,1+\alpha}([0,T]\times\bar D)$ be uniformly with
respect to $u\in[\sigma,\omega]$ and $f(x,0,0)=0$ on $\partial D$.
Fixed $p>\max\{d,1+\frac{d}{2}\}$. If there exists at least one $u_0\in W_{p,\mathscr{B}}^{2}$ satisfying
$\underline u(0,x)\leq u_0(x)\leq\bar u(0,x)$, then the problem (\ref{3.8}) has at least one solution
$u\in C^{1+\alpha/2,2+\alpha}([0,T]\times\bar D)$ and satisfies
   \bess
\underline u\leq u\leq\bar u\ \ \ \ on\ \ \ [0,T]\times\bar D.
  \eess
   \end{prop}
The proof of Proposition \ref{p3.1} is standard and we omit it here.
Now, by using the Proposition \ref{p3.1}, we only need find a pair of upper
and lower solution to the problem (\ref{3.8}). In order to do that, we need consider
the periodic parabolic eigenvalue problem
       \bes\left\{\begin{array}{llll}
  \partial_tu(t,x)+\mathcal{A}(t)u(t,x)=\lambda u, \ \ &in\  D\times(0,T],\\
\mathscr{B}u=0, \ \ \ \ &on\ \partial D\times(0,T],\\
u(0,x)=u(T,x), \ \ \ & in \ D.
   \end{array}\right.\lbl{3.9}\ees
It is easy to check that if $b\geq0$ on $\partial D\times(0,T]$, then the
maximum principle holds for the problem $(\mathcal{A},\mathscr{B})$.
We will need the following Lemmas.
  \begin{lem}\lbl{l3.2} \cite[Proposition14.4]{Hess1991book} The  principal
eigenvalue $\lambda_1$ of (\ref{3.9}) exists uniquely. Furthermore, if
$a_0\geq0$, and $a_0\not\equiv0$ when $\mathscr{B}=\partial_\nu$, then $\lambda_1>0$. In case $a_0=0$ and $\mathscr{B}=\partial_\nu$, we have
$\lambda_1=0$.
   \end{lem}
We want to know how the principle eigenvalue $\lambda_1$ depends on the
zero-order term $a_0$. Because in Fokker-Planck equations, the role of drift
term is reflected in the zero-order term $a_0$. In order to do that, we need the following lemma.
  \begin{lem}\lbl{l3.3} \cite[Proposition16.6]{Hess1991book} For the
inhomogeneous problem
  \bes
\partial_tu(t,x)+\mathcal{A}(t)u(t,x)-\lambda u=h, \ \ h\geq,\not\equiv0,
   \lbl{3.100}\ees
where $h\in\mathbb{F}$ (see Section 2 for the definition of $\mathbb{F}$).
Let $\lambda_1$ be the  principal
eigenvalue $\lambda_1$ of (\ref{3.9}). Then we have

(i) If $\lambda<\lambda_1$, then the problem (\ref{3.100}) has a unique
solution $u$ and $u>0$ in $\mathbb{F}_1$;

(ii) If $\lambda\geq\lambda_1$, then the problem (\ref{3.100}) has
no positive solution, and no solution at all if $\lambda=\lambda_1$.
   \end{lem}

By using the above Lemma \ref{l3.3}, it is easy to prove the following Lemma.
  \begin{lem}\lbl{l3.4} Let $\lambda_1=\lambda_1(a_0)$ be the  principal
eigenvalue of (\ref{3.9}). Then $\lambda(a_0)$ is strictly increasing
in $a_0$.
   \end{lem}

{\bf Proof.} suppose $a_0^1(t,x)\leq a_0^2(t,x)$ and
$a_0^1(t,x)\not\equiv a_0^2(t,x)$. Suppose $\phi_i$ is the corresponding
positive eigenfunction to $\lambda_1(a^i_0)$ and $\lambda_1(a^i_0)>0$,
$i=1,2$. Denote $\mathcal{A}_0=\mathcal{A}-a_0$.
We aim to prove that $\lambda_1(a^1_0)<\lambda_1(a^2_0)$.
 On the
contrary, we suppose $\lambda_1(a^1_0)\geq\lambda_1(a^2_0)$, then we have
   \bess
&&\partial_t(\phi_1-\phi_2)+\mathcal{A}_0(t)(\phi_1-\phi_2)+a_0^2(\phi_1-\phi_2)\\
&=&\partial_t\phi_1+\mathcal{A}_0(t)\phi_1+a_0^2\phi_1-\lambda_1(a^2_0)\phi_2\\
&\geq&\partial_t\phi_1+\mathcal{A}_0(t)\phi_1+a_0^1\phi_1-\lambda_1(a^2_0)\phi_2\\
&=&\lambda_1(a^1_0)\phi_1-\lambda_1(a^2_0)\phi_2\\
&\geq&\lambda_1(a^1_0)(\phi_1-\phi_2),
   \eess
that is to say,
   \bess
\partial_t(\phi_1-\phi_2)+\mathcal{A}_0(t)(\phi_1-\phi_2)+a_0^2(\phi_1-\phi_2)
-\lambda_1(a^1_0)(\phi_1-\phi_2)
=:h\geq,\not\equiv0.
  \eess
By using comparison principle, we deduce that $\phi_1-\phi_2>0$
 if $\lambda_1(a_0^1)>\lambda_1(a_0^2)$. We obtain
 a contradiction with (ii) of Lemma \ref{l3.3}.  $\Box$

Now, we use the above discussion to solve the problem (\ref{3.6}).
  \begin{theo}\lbl{t3.2} Suppose there exist two positive constants $c_0$ and $M_0$ satisfying
   \bess
a_0(t,x)\geq c_0, \ \ f(t,x,0)>0,\ \ f(t,x,\xi)\leq0\ {\rm for}\ \xi\geq M_0,\ \  f\in\mathbb{F}.
   \eess
Then the problem (\ref{3.6}) admits a unique solution.
  \end{theo}

 {\bf Proof.} The existence of periodic solution is obtained by using
 Proposition \ref{p3.1}. We only need find a pair of upper-lower solution of (\ref{3.6}).
Actually, from the assumptions on $a_0$ and $f$, we see that $\bar u=M\geq M_0$ is an upper
solution. Let $\phi$ be a positive eigenfunction corresponding to $\lambda_1(a_0)$,
i.e., $\phi$ is the solution of (\ref{3.9}) with $\lambda=\lambda_1(a_0)$.
Take $\varepsilon>0$ and set $\underline u=\varepsilon\phi$, then $\underline u$ is s lower
solution of (\ref{3.6}). Moreover, $\underline u$ and $\bar u$ are the ordered upper and
lower solutions of (\ref{3.6}) if we choose $\varepsilon\ll1$ and $M\gg1$.
According to Proposition \ref{p3.1}, the problem (\ref{3.6}) has at least one solution
$u$ satisfying $\varepsilon\phi\leq u\leq M$. The proof of uniqueness follows from the
comparison principle. $\Box$

\begin{remark}\lbl{r3.2} Following Lemma \ref{l3.2}, we can see that if $a_0=0$, then the problem (\ref{3.9})
with Neumann boundary admits the principle eigenvalue $\lambda_1=0$ and eigenfunction $\phi=1$.

A typical example in \ref{t3.2} is $f(t,x)=f_1(t,x)-uf_2(t,x)$ with $f_i(t,x)\in(N_1,N_2)$ for all
$(t,x)\in[0,T]\times \bar D$ and $N_1>0$.

The assumptions on $f$ can be given weaker, but it is not our aim. See \cite{Wang2016} for
$f=u(h_1(t,x)-h_2(t,x)u)$ and $h_i$, $i=1,2$, are some functions.
\end{remark}
\subsection{Bounded domain with reflecting boundary condition}
It is easy to see that there is no periodic solution to equation (\ref{3.2}). Now, we
consider another case.
We assume that a stochastic process $X_t$ satisfies
   \bes\left\{\begin{array}{llll}
dX_t=b(t,X_t)dt+\sigma(t,X_t)dB_t,\\
 X_0=x,
   \end{array}\right.\lbl{3.10}\ees
with the reflecting boundary condition \cite[Section 5.1.1]{Gardiner2009book},
where $b$ is an $d$-dimensional vector function, $\sigma$ is an $d\times m$ matrix function and
$B_t$ is an $m$-dimensional Brownian motion.
Our aim is to study the  existence of periodic solution to the
corresponding Fokker-Planck equation, i.e.,  the existence of the following equations
   \bes\left\{\begin{array}{llll}
  \partial_tp(t,x)=\mathcal{A}^*p(t,x), \ \ &in\  D\times(0,T],\\
(b\cdot\nu) p-p_{\partial D}\cdot\nu=0, \ \ \ \ &on\ \partial D\times(0,T],\\
p(0,x)=p(T,x), \ \ \ & in \ D,
   \end{array}\right.\lbl{3.11}\ees
where
   \bess
\mathcal{A}^*p&=&-\sum_i\frac{\partial}{\partial x_i}(b_ip)+\frac{1}{2}\sum_{i,j}
\frac{\partial^2}{\partial x_i\partial x_j}((\sigma\sigma^T)_{ij}p),\\
p_{\partial D}&=&(\sum_i\frac{\partial}{\partial x_i}((\sigma\sigma^T)_{i1}p),\cdots,\sum_i\frac{\partial}{\partial x_i}((\sigma\sigma^T)_{id}p)).
  \eess
It follows the results of \cite{Duan2015} that for the existence of probability density function,
the necessary condition is that operator $\mathcal{A}^*$ is a uniform elliptic operator.
Throughout this section we assume that the operator $\mathcal{A}^*$ is a uniform elliptic operator.

We first consider a special case. It is noted that most of work on the
periodic parabolic problem the authors assumed the boundary function $b$
does not depend on the time $t$. Because under this assumption, one can
apply the standard theory of evolution equation of "parabolic type", see
\cite{Amann1978,Wang2016}. So we assume
   \bess
(\sigma\sigma^T)_{ij}(t,x)=\alpha(t),\ \ 0<\alpha_0\leq\alpha(t)\leq\alpha_1,\ \
b(t,x)=\alpha(t)\times b_0(x).
  \eess
Then the problem (\ref{3.11}) becomes
   \bes\left\{\begin{array}{llll}
  \partial_tp(t,x)=\mathcal{A}^*p(t,x), \ \ &in\  D\times(0,T],\\
\partial_\nu p-b_0p=0, \ \ \ \ &on\ \partial D\times(0,T],\\
p(0,x)=p(T,x), \ \ \ & in \ D.
   \end{array}\right.\lbl{3.12}\ees
We can use the similar method to deal with the problem (\ref{3.12}).

It is well known that the upper-lower method is not suitable
to the linear parabolic equation. The reason is that if we find an upper solution $\phi$ for
a linear parabolic equation, then $\lambda\phi$ will be an upper solution for any $\lambda>0$.
Hence we can not obtain the existence of non-negative solution for this linear parabolic equation.
And in this subsection, we only consider the one dimensional case because it can be
calculated clearly. We want to obtain the existence of periodic solution of (\ref{3.11}).
For simplicity, we denote $a(t,x)=(\sigma\sigma^T)(t,x)$ and $D=(0,1)$.
Due to the operator $\mathcal{A}^*$ is a uniform elliptic operator, we have $a(t,x)>0$ for
$(t,x)\in[0,T]\times[0,1]$.
The one dimensional problem will be written as
    \bes\left\{\begin{array}{llll}
  \partial_tp(t,x)-(a(t,x)p(t,x))_{xx}+(b(t,x)p(t,x))_x=0, \ \ &in\  D\times(0,T],\\
(a(t,x)p(t,x))_x -bp=0, \ \ \ \ &on\ \partial D\times(0,T],\\
p(0,x)=p(T,x), \ \ \ & in \ D.
   \end{array}\right.\lbl{3.13}\ees
We first assume that
   \bess
(a(t,x)p(t,x))_x -bp=0, \ \ \ {\rm in}\ \ [0,T]\times[0,1].
  \eess
Then we get
  \bes
p(t,x)=\exp\left(\int\frac{b-a_x}{a}dx\right),
   \lbl{3.14}\ees
which implies that $p_t=0$, i.e.,
   \bes
\int\frac{a(b_t-a_{xt})-a_t(b-a_x)}{a^2}dx=0.
  \lbl{3.15}\ees
That is to say, the stochastic process has stationary probability measure.
Summing the above discussion, we have
   \begin{theo}\lbl{t3.3} Suppose (\ref{3.15}) hold. Then problem (\ref{3.13}) admits
a solution $p$ satisfying (\ref{3.14}).
  \end{theo}

For $d\geq2$, we can not calculate it clearly. But we guess there
exists a positive periodic solution to problem (\ref{3.11}). Indeed,
it follows from (\ref{3.11}) that
   \bes
\int_Dp_0(x)dx=\int_Dp(t,x)dx, \ \ \ \forall t>0.
  \lbl{3.16}\ees
The existence of periodic solution to (\ref{3.11}) is equivalent to getting
$p(0,x)=p(T,x)$ point by point for $x\in D$ from (\ref{3.16}).

In 2000, Lieberman \cite{Lieberman2001jmaa,Lieberman2000,Lieberman2001na} did
a series of work about the periodic solution of parabolic equation on bounded
domain. Especially in \cite{Lieberman2001na}, Lieberman obtained the existence of
periodic of the following parabolic equation
  \bess
\left\{\begin{array}{llll}
u_t-div A(t,x,u,\nabla u)+B(t,x,u,\nabla u)=0, \ in \ (0,T)\times\partial D,\\
A(t,x,u,\nabla u)\cdot\nu+\psi(t,x,u)=0,\ \ \ on \ (0,T)\times\partial D,\\
u(0,x)=u(T,x),\ \  in\ D.
  \end{array}\right.
  \eess
If $b,\,\sigma$ satisfy the conditions of \cite[Lemma 2.1]{Lieberman2001na}, then
the (\ref{3.11}) will admit a periodic solution $p$.

\subsection{Whole space}
In this subsection, we consider the existence of periodic solutions of Fokker-Planck equations in the whole space.
For a stochastic process $X_t$ satisfies equation (\ref{3.1}), the corresponding Fokker-Planck equation is the following form
   \bes
\partial_tp=-\sum_i\frac{\partial}{\partial x_i}(b_ip)+\frac{1}{2}\sum_{i,j}
\frac{\partial^2}{\partial x_i\partial x_j}((\sigma\sigma^T)_{ij}p).
  \lbl{3.17}\ees
Furthermore, if the probability density $p(t,x)$ satisfies $p(t+T,x)=p(t,x)$, $\forall(t,x)$, then
$p(t,x)$ is call a $T$-periodic solution of (\ref{3.17}).

By using the method of \cite{Georgiev2013}, we will obtain the existence of periodic solution
of (\ref{3.17}). In \cite{Georgiev2013}, the author considered the following periodicity problem
   \bes\left\{\begin{array}{llll}
u_t-\Delta u=f(t,x,u,u_x), \ \ \ & t>0,\ \ x\in\mathbb{R}^d,\\
u(t,x)=u(t+T,x), \ \ \ & t\geq0,\ \ x\in\mathbb{R}^d,\\
  \end{array}\right.\lbl{3.18}\ees
where $d\geq2$, $f\in C(\mathbb{R},\mathbb{R}^d,\mathbb{R},\mathbb{R}^d)$, $u_x=(u_{x_1},u_{x_2},\cdots,u_{x_d})$,
$f$ is $T$-periodic function with respect to the time variable $t$, the period $T>0$ is arbitrary chosen and fixed.
They got the following result.
 \begin{prop}\lbl{p3.2} Let $d\geq2$, $n\in\mathbb{N}$ be fixed, $T>0$ be fixed,
$f\in C(\mathbb{R},\mathbb{R}^d,\mathbb{R},\mathbb{R}^d)$. $f$ is $T$-periodic with respect to the time variable
$t$. Also let $0\leq c_i,l_i,m_i,p_i,q_i,l_i<\infty$, $i=1,2,\cdots,n$, be fixed constants, $0\leq k_i,n_i<\infty$,
$i=1,2,\cdot,d$, be fixed constants, $b_i(t)\in C(\mathbb{R}_+),\ g_i(x)\in C(\mathbb{R}^d)$,
$\sup_{\mathbb{R}_+}|b_i(t)|<\infty$, $\sup_{\mathbb{R}^d}|g_i(x)|<\infty$, $i=1,2,\cdots,n$,
  \bess
|f(t,x,u,u_x)|\leq\sum_{i=1}^n\left(c_i|b_i(t)|^{p_i}+l_i|u|^{q_i}+m_i|g_i(x)|^{l_i}\right)+\sum_{i=1}^dk_i|u_{x_i}|^{n_i}
   \eess
for every $(t,x,u,u_x)\in(\mathbb{R},\mathbb{R}^d,\mathbb{R},\mathbb{R}^d)$. Then the problem (\ref{3.18}) has a solution
$u\in\mathcal{C}^1(\mathbb{R}_+,\mathcal{C}^2(\mathbb{R}^d))$ ($\mathcal{C}^1(\mathbb{R}_+,\mathcal{C}^2(\mathbb{R}^d))$
will be defined later).
\end{prop}

Comparing the problem (\ref{3.17}) and (\ref{3.18}), we see that the problem (\ref{3.17}) is
a linear problem and problem (\ref{3.18}) will contain problem (\ref{3.17}) if $(\sigma\sigma^T)_{ij}=constant$.
Firstly, it is remarked that when $d=1$, the results in subsection 3.2 also holds for problem
(\ref{3.17}). That is to say, theorem \ref{t3.3} holds for (\ref{3.17}). In order to get
the existence of solutions to (\ref{3.17}) for $d\geq2$, we suppose that $a(t)$ is
continuous positive $T$-periodic function, which is defined on the whole real axis $\mathbb{R}$.
We denote  $[a]=\frac{1}{T}\int_0^Ta(t)dt$. For $D\subset\mathbb{R}^d$, set
   \bess
\mathcal{C}^1([0,T],\mathcal{C}^2(D))&=&\{u(t,x):\ continuously-differentiable \ in \ t\in\mathbb{R}_+,\\
&&\qquad\qquad twice\
continuously-differentiable\ in\ x\in D, \\
&&\qquad\qquad  u(t+T,x)=u(t,x)\ for\ t\geq0\ and\ x\in D\}.
   \eess
Fix $0<Q<\infty,\ 0<\varepsilon<1$ and denote
   \bess
F&=&\max\left\{\sup_{t\geq0,x\in\mathbb{R}^d}div b(t,x),\ \sup_{t\geq0,x\in\mathbb{R}^d}\sum_{i=1}^d|b_i(t,x)|, \ \
\sup_{t\geq0,x\in\mathbb{R}^d}\sum_{i,j=1}^d|(\sigma\sigma^T)_{ij}(t,x)|\right\},\\
G&=&\max\left\{\sup_{0\leq t,s\leq T}\frac{e^{-[a]T}}{1-e^{-[a]T}}e^{\int_t^{t+s}a(r)dr},\ \sup_{0\leq t\leq T}a(t), \ \
\sup_{0\leq t\leq T}e^{\int_0^ta(r)dr}\right\},\\
S&=&(2d+1)(F^2+Q^2).
  \eess
If we assume $\sup_{t,x}|b_i(t,x)|<\infty$,
$\sup_{t,x}|(\sigma\sigma^T)_{ij}(t,x)|<\infty$, $i,j=1,2,\cdots,d$, we have that
$0\leq F<\infty$. Note that $a(t)$ is a continuous positive $T$-periodic function,
which is defined on the whole real axis $\mathbb{R}$, we conclude that $0\leq G<\infty$.
From here and $0<Q<\infty$ we get $0\leq S<\infty$.

The main result of this subsection is the following theorem.
  \begin{theo}\lbl{t3.4} Assume that the SDEs (\ref{3.1})
admits a probability density function $p(t,x)$ satisfying (\ref{3.17}). Assume further
that $F<\infty$ and $b_i,\ \sigma$ are $T$-periodic functions, then there exists a periodic solution
$p(t,x)\in\mathcal{C}^1(\mathbb{R}_+,\mathcal{C}^2(\mathbb{R}^d))$
satisfying $(\ref{3.17})$ with $p(t+T,x)=p(t,x)$ for $t\geq0$ and $x\in\mathbb{R}^d$.
 \end{theo}

We choose the constants $A_i$, $i=1,2,\cdots,d$, so that $A_i>0$ satisfying
\bes
&&(A_1\cdots A_d)^2Q+GT\Big[G(A_1\cdots A_d)^2Q+\sum_{i=1}^d(A_1\cdots A_{i-1}A_{i+1}\cdots A_d)^2SQ\nonumber\\
&&+
\sum_{i,j=1,i\neq j}^dA_iA_j(A_1\cdots A_{i-1}A_{i+1}\cdots A_{j-1}A_{j+1}\cdots A_d)^2SQ\nonumber\\
&&+
\sum_{i=1}^dA_i(A_1\cdots A_{i-1}A_{i+1}\cdots A_d)^2SQ\Big]\leq(1-\varepsilon)Q,\nonumber\\
&&(A_1\cdots A_d)^2Q+(G^3T+1)\Big[G(A_1\cdots A_d)^2Q+\sum_{i=1}^d(A_1\cdots A_{i-1}A_{i+1}\cdots A_d)^2SQ\nonumber\\
&&+
\sum_{i,j=1,i\neq j}^dA_iA_j(A_1\cdots A_{i-1}A_{i+1}\cdots A_{j-1}A_{j+1}\cdots A_d)^2SQ\nonumber\\
&&+
\sum_{i=1}^dA_i(A_1\cdots A_{i-1}A_{i+1}\cdots A_d)^2SQ\Big]\leq(1-\varepsilon)Q,\lbl{3.19}\\
&&(A_1\cdots A_{k-1}A_{k+1}\cdots A_d)^2A_kQ+(G^3T+1)\Big[GT(A_1\cdots A_{k-1}A_{k+1}\cdots A_d)^2A_kQ\nonumber\\
&&+\sum_{i=1,i\neq k}^d(A_1\cdots A_{i-1}A_{i+1}\cdots A_{k-1}A_{k+1}\cdots A_d)^2A_kSQ+(A_1\cdots A_{i-1}A_{i+1}\cdots A_d)^2SQ\nonumber\\&&+
\sum_{i,j=1,i\neq j\neq k}^dA_iA_jA_k(A_1\cdots A_{i-1}A_{i+1}\cdots A_{j-1}A_{j+1}\cdots A_{k-1}A_{k+1}\cdots A_d)^2SQ\nonumber\\
&&+2\sum_{j,k=1,k\neq j}^dA_kA_j(A_1\cdots A_{j-1}A_{j+1}\cdots A_{k-1}A_{k+1}\cdots A_d)^2SQ\nonumber\\
&&+
\sum_{i=1,i\neq k}^dA_iA_k(A_1\cdots A_{i-1}A_{i+1}\cdots A_{k-1}A_{k+1}\cdots A_d)^2SQ\nonumber\\
&&+(A_1\cdots A_{i-1}A_{i+1}\cdots A_d)^2SQ\Big]\leq(1-\varepsilon)Q \ \ \ \forall j=1,2,\cdots,d,\nonumber\\
&&A_0=A_{n+1}=1.\nonumber
  \ees
Such choice is suitable if $0<A_i<1$, $i=1,2,\cdots,d$, is small enough. We set
$A=(A_1,A_2,\cdots,A_d)$ and
  \bess
B_1=\{x\in \mathbb{R}^d:\ 0\leq x_i\leq A_i, \ i=1,2,\cdots,d\}.
  \eess
We first prove that the periodicity problem
   \bes\left\{\begin{array}{llll}
\partial_tp=-\sum_i\frac{\partial}{\partial x_i}(b_ip)+\frac{1}{2}\sum_{i,j}
\frac{\partial^2}{\partial x_i\partial x_j}((\sigma\sigma^T)_{ij}p), \ \ \ & t>0,\ \ x\in B_1,\\
p(t,x)=p(t+T,x), \ \ \ & t\geq0,\ \ x\in B_1.
  \end{array}\right.\lbl{3.20}\ees
We will use the fixed point arguments to prove the existence of solution to (\ref{3.20}).
Following the idea of \cite{Georgiev2013}, we give the following lemma.
   \begin{lem}\lbl{l3.5} If $p\in \mathcal{C}^1([0,T],\mathcal{C}^2(B_1))$ satisfies the equation
      \bes
&&\int_x^A\int_y^Ap(t,z)dzdy-\frac{e^{-[a]T}}{1-e^{-[a]T}}\int_0^Te^{\int_t^{t+s}a(r)dr}\Big(a(t+s)
\int_x^A\int_y^A\nonumber\\
&& p(t+s,z)dzdy+\sum_{i=1}^d\int_{\bar x_i}^A\int_{\bar y_i}^A
(\sigma\sigma^T)_{ii}(t+s,\hat z_i)p(t+s,\hat z_i)d\hat z_id\hat y_i\nonumber\\
&&+\sum_{i,j=1,i\neq j}^d\int_{x}^A\int_{\bar y_{ij}}^A
(\sigma\sigma^T)_{ij}(t+s,\hat z_{ij})p(t+s,\hat z_{ij})d\hat z_{ij}dy\nonumber\\
&&+\sum_{i=1}^d\int_x^A\int_{\breve y_i}^A[b_i(t+s,\breve z_i)p(t+s,\breve z_i)]d\breve z_idy\Big)ds=0,
   \lbl{3.21}\ees
then $p(t,x)$ is a solution to the problem (\ref{3.20}). We use the following symbols
  \bess
&&\int_x^A=\int_{x_1}^{A_1}\int_{x_2}^{A_2}\cdots\int_{x_d}^{A_d}, \ \ \ \
\int_{\bar x_i}^A=\int_{x_1}^{A_1}\cdots\int_{x_{i-1}}^{A_{i-1}}\int_{x_{i+1}}^{A_{i+1}}\cdots\int_{x_d}^{A_d},\\
&&\int_{\bar y_{ij}}^A=\int_{y_1}^{A_1}\cdots\int_{y_{i-1}}^{A_{i-1}}\int_{y_{i+1}}^{A_{i+1}}\cdots
\int_{y_{j-1}}^{A_{j-1}}\int_{y_{j+1}}^{A_{j+1}}\cdots\int_{x_d}^{A_d},\\
&&\int_{\breve y_i}^A=\int_{y_1}^{A_1}\cdots\int_{y_{i-1}}^{A_{i-1}}\int_{y_{i+1}}^{A_{i+1}}\cdots\int_{y_d}^{A_d},
\  \ \breve z_i=(z_1,\cdots,z_{i-1},y_i,z_{i+1},z_d),\\
&&\hat w_i=(w_1,\cdots,w_{i-1},x_i,w_{i+1},\cdots,w_d), \ \ d\hat w_i=dw_1\cdots dw_{i-1}dw_{i+1}\cdots dw_n,\ \ w=y\ {\rm or}\ z,\\
&& \hat z_{ij}=(z_1,\cdots,z_{i-1},y_i,z_{i+1},\cdots,z_{j-1},y_j,z_{j+1},z_d).
   \eess
   \end{lem}

{\bf Proof.} Differentiating the (\ref{3.21}) twice in $x_1$, twice in $x_2,\cdots$, twice in $x_d$ and using the
periodicity of $a,p,b$ and $\sigma$, we can obtain the desired result. See \cite[Lemma 2.1]{Georgiev2013} for more details.
 $\Box$

Lemma \ref{l3.5} implies that the existence of solution to (\ref{3.20}) is equivalent to
the existence of fixed point of $L_1$, i.e., $L_1(p)=p$, where
  \bess
L_1(p)&=&p(t,x)+\int_x^A\int_y^Ap(t,z)dzdy-\frac{e^{-[a]T}}{1-e^{-[a]T}}\int_0^Te^{\int_t^{t+s}a(r)dr}\Big(a(t+s)
\int_x^A\int_y^A\nonumber\\
&& p(t+s,z)dzdy+\sum_{i=1}^d\int_{\bar x_i}^A\int_{\bar y_i}^A
(\sigma\sigma^T)_{ii}(t+s,\hat z_i)p(t+s,\hat z_i)d\hat z_id\hat y_i\nonumber\\
&&+\sum_{i,j=1,i\neq j}^d\int_{x}^A\int_{\bar y_{ij}}^A
(\sigma\sigma^T)_{ij}(t+s,\hat z_{ij})p(t+s,\hat z_{ij})d\hat z_{ij}dy\nonumber\\
&&+\sum_{i=1}^d\int_x^A\int_{\breve y_i}^A[b_i(t+s,\breve z_i)p(t+s,\breve z_i)]d\breve z_idy\Big)ds.
  \eess
In order to get the fixed point of $L_1$, we define
   \bess
D_1&=&\{u\in\mathcal{C}^1([0,T],\mathcal{C}^2(B_1)):\ |u|\leq Q,|u_t|\leq Q,|u_{x_i}|\leq Q, \ i=1,2,\cdots,d\},\\
\tilde D_1&=&\{u\in\mathcal{C}^1([0,T],\mathcal{C}^2(B_1)):\ |u|\leq(1+\varepsilon) Q,|u_t|\leq (1+\varepsilon)Q,\\
&&\qquad|u_{x_i}|\leq (1+\varepsilon)Q, \ i=1,2,\cdots,d\}.
  \eess
In the set $D_1$ and $\tilde D_1$, we define a norm as follows:
  \bess
\|u\|=\max\left\{\max_{t\in[0,T],x\in B_1}|u|,\ \max_{t\in[0,T],x\in B_1}|u_t|,\ \max_{t\in[0,T],x\in B_1}|u_{x_i}|,\ i=1,2,\cdots,d\right\}.
  \eess
Then $D_1$, $\tilde D_1$ and $\mathcal{C}^1([0,T],\mathcal{C}^2(B_1))$ are completely normed spaces with
respect to this norm, see Appendix of \cite{Georgiev2013}. We rewrite the operator $L_1$ in the following form
  \bess
L_1(p)=M_1(p)+N_1(p),
  \eess
where
  \bess
M_1(p)&=&(1+\varepsilon)p,\\
N_1(p)&=&-\varepsilon p+\int_x^A\int_y^Ap(t,z)dzdy-\frac{e^{-[a]T}}{1-e^{-[a]T}}\int_0^Te^{\int_t^{t+s}a(r)dr}\Big(a(t+s)
\int_x^A\int_y^A\nonumber\\
&& p(t+s,z)dzdy+\sum_{i=1}^d\int_{\bar x_i}^A\int_{\bar y_i}^A
(\sigma\sigma^T)_{ii}(t+s,\hat z_i)p(t+s,\hat z_i)d\hat z_id\hat y_i\nonumber\\
&&+\sum_{i,j=1,i\neq j}^d\int_{x}^A\int_{\bar y_{ij}}^A
(\sigma\sigma^T)_{ij}(t+s,\hat z_{ij})p(t+s,\hat z_{ij})d\hat z_{ij}dy\nonumber\\
&&+\sum_{i=1}^d\int_x^A\int_{\breve y_i}^A[b_i(t+s,\breve z_i)p(t+s,\breve z_i)]d\breve z_idy\Big)ds.
  \eess
To obtain the operator $L_1$ has a fixed point in the space $\mathcal{C}^1([0,T],\mathcal{C}^2(B_1))$ we
need the following lemma.
  \begin{lem}\lbl{l3.6}\cite[Corollary 2.4, p.3231]{XiangYuan2009} Let $X$ be a nonempty closed
convex subset of a Banach space $Y$. Suppose that $T$ and $S$ map $X$ into $Y$ such that

(i) $S$ is continuous, $S(X)$ resides in a compact subset of $Y$;

(ii) $T:X\to Y$ is expansive and onto.

Then there exists a point $x^*\in X$ with $Sx^*+Tx^*=x^*$.
  \end{lem}

We recall the definition of expansive operator.
 \begin{defi}\lbl{d3.2} \cite{XiangYuan2009} Let $(X,d)$ be a metric space and
$M$ be a subset of $X$. The mapping $T:M\to X$ is said to be expansive, if there
exists a constant $h>1$ such that
  \bess
d(Tx,Ty)\geq hd(x,y), \ \ \ \ \forall x,y\in M.
  \eess
 \end{defi}
It is easy to check the following lemma, see \cite[lemma 2.3]{Georgiev2013} for the details.
  \begin{lem}\lbl{l3.7} The operator $M_1:D_1\to\tilde D_1$ is an expansive operator and onto.
    \end{lem}
Next we prove the operator $N_1$ satisfies the (i) of Lemma \ref{l3.6}.
  \begin{lem}\lbl{l3.8} The operator $N_1:D_1\to D_1$ is continuous and $D_1$ is a compact set in $\tilde D_1$.
    \end{lem}

{\bf Proof.} We first prove $N_1$ maps $D_1$ to $D_1$. For any $p\in D_1$, by using (\ref{3.19}), we have
   \bess
|N_1(p)|&\leq&\varepsilon |p|+\int_x^A\int_y^A|p(t,z)|dzdy+\frac{e^{-[a]T}}{1-e^{-[a]T}}\int_0^Te^{\int_t^{t+s}a(r)dr}\Big(a(t+s)
\int_x^A\int_y^A\nonumber\\
&& |p(t+s,z)|dzdy+\sum_{i=1}^d\int_{\bar x_i}^A\int_{\bar y_i}^A
|(\sigma\sigma^T)_{ii}(t+s,\hat z_i)p(t+s,\hat z_i)|d\hat z_id\hat y_i\nonumber\\
&&+\sum_{i,j=1,i\neq j}^d\int_{x}^A\int_{\bar y_{ij}}^A
|(\sigma\sigma^T)_{ij}(t+s,\hat z_{ij})p(t+s,\hat z_{ij})|d\hat z_{ij}dy\nonumber\\
&&+\sum_{i=1}^d\int_x^A\int_{\breve y_i}^A|b_i(t+s,\breve z_i)p(t+s,\breve z_i)|d\breve z_idy\Big)ds\\
&\leq&\varepsilon Q+(A_1\cdots A_d)^2Q+GT\Big[G(A_1\cdots A_d)^2Q\\
&&+\sum_{i=1}^d(A_1\cdots A_{i-1}A_{i+1}\cdots A_d)^2SQ\\
&&+
\sum_{i,j=1,i\neq j}^dA_iA_j(A_1\cdots A_{i-1}A_{i+1}\cdots A_{j-1}A_{j+1}\cdots A_d)^2SQ\\
&&+
\sum_{i=1}^dA_i(A_1\cdots A_{i-1}A_{i+1}\cdots A_d)^2SQ\Big]\\
&\leq&\varepsilon Q+(1-\varepsilon)Q
  \eess
for every $t\in[0,T]$ and every $x\in B_1$.

For $(N_1(p))_t$, we get
   \bess
(N_1(p))_t&=&-\varepsilon p_t+\int_x^A\int_y^Ap_t(t,z)dzdy+\frac{e^{-[a]T}}{1-e^{-[a]T}}
a(t)e^{-\int_0^ta(r)dr}\int_t^{t+T}e^{\int_0^{t_1}a(r)dr}\nonumber\\
&&\times\Big[a(t_1)
\int_x^A\int_y^A p(t_1,z)dzdy+\sum_{i=1}^d\int_{\bar x_i}^A\int_{\bar y_i}^A
(\sigma\sigma^T)_{ii}(t_1,\hat z_i)p(t_1,\hat z_i)d\hat z_id\hat y_i\nonumber\\
&&+\sum_{i,j=1,i\neq j}^d\int_{x}^A\int_{\bar y_{ij}}^A
(\sigma\sigma^T)_{ij}(t_1,\hat z_{ij})p(t_1,\hat z_{ij})d\hat z_{ij}dy\nonumber\\
&&+\sum_{i=1}^d\int_x^A\int_{\breve y_i}^A[b_i(t_1,\breve z_i)p(t_1,\breve z_i)]d\breve z_idy\Big]dt_1\\
&&+\Big[a(t)
\int_x^A\int_y^Ap(t,z)dzdy+\sum_{i=1}^d\int_{\bar x_i}^A\int_{\bar y_i}^A
(\sigma\sigma^T)_{ii}(t,\hat z_i)p(t,\hat z_i)d\hat z_id\hat y_i\nonumber\\
&&+\sum_{i,j=1,i\neq j}^d\int_{x}^A\int_{\bar y_{ij}}^A
(\sigma\sigma^T)_{ij}(t,\hat z_{ij})p(t,\hat z_{ij})d\hat z_{ij}dy\nonumber\\
&&+\sum_{i=1}^d\int_x^A\int_{\breve y_i}^A[b_i(t,\breve z_i)p(t,\breve z_i)]d\breve z_idy\Big]
  \eess
which implies that (using (\ref{3.19}) again)
   \bess
|(N_1(p))_t|&\leq&
  \varepsilon Q+(A_1\cdots A_d)^2Q+(G^3T+1)\Big[G(A_1\cdots A_d)^2Q\\
&&+\sum_{i=1}^d(A_1\cdots A_{i-1}A_{i+1}\cdots A_d)^2SQ\\
&&+
\sum_{i,j=1,i\neq j}^dA_iA_j(A_1\cdots A_{i-1}A_{i+1}\cdots A_{j-1}A_{j+1}\cdots A_d)^2SQ\\
&&+
\sum_{i=1}^dA_i(A_1\cdots A_{i-1}A_{i+1}\cdots A_d)^2SQ\Big]\\
&\leq&\varepsilon Q+(1-\varepsilon)Q, \ \ \forall t\in[0,T],\ x\in\mathbb{R}^d.
  \eess
Let $k\in\{1,2,\cdots,d\}$ be arbitrary chosen and fixed. Then we have
  \bess
(N_1(p))_{x_k}&=&-\varepsilon p_{x_k}-\int_{\bar x_k}^A\int_{\hat y_k}^Ap_t(t,z)dzd\hat y_k\\
&&
-\frac{e^{-[a]T}}{1-e^{-[a]T}}\int_0^Te^{\int_t^{t+s}a(r)dr}\Big(-a(t+s)
\int_{\bar x_k}^A\int_{\hat y_k}^Ap(t+s,z)dzd\hat y_k\nonumber\\
&& +\sum_{i=1,i\neq k}^d\int_{\overline{(\bar x_i)_k}}^A\int_{\widehat{(\bar y_i)_k}}^A
(\sigma\sigma^T)_{ii}(t+s,\hat z_i)p(t+s,\hat z_i)d\hat z_id\widehat{(\bar y_i)_k} \nonumber\\
&&+\int_{\bar x_k}^A\int_{\bar y_k}^A
(\sigma\sigma^T)_{kk}(t+s,\hat z_k)p(t+s,\hat z_k)d\hat z_kd\hat y_k\\
&&+\sum_{i,j=1,i\neq j\neq k}^d\int_{\bar x_k}^A\int_{\widehat{(\bar y_{ij})_k}}^A
(\sigma\sigma^T)_{ij}(t+s,\hat z_{ij})p(t+s,\hat z_{ij})d\hat z_{ij}d\hat y_k\nonumber\\
&&+\sum_{i,j=1,i\neq j,j=k}^d\int_{\bar x_k}^A\int_{\bar y_{ik}}^A
(\sigma\sigma^T)_{ik}(t+s,\hat z_{ik})p(t+s,\hat z_{ik})d\hat z_{ik}d\hat y_k\nonumber\\
&&+\sum_{i,j=1,i\neq j,i=k}^d\int_{\bar x_k}^A\int_{\bar y_{kj}}^A
(\sigma\sigma^T)_{kj}(t+s,\hat z_{kj})p(t+s,\hat z_{kj})d\hat z_{kj}d\hat y_k\nonumber\\
&&+\sum_{i=1,i\neq k}^d\int_{x_k}^A\int_{\widehat{(\breve y_i)_k}}^A[b_i(t+s,\widehat{(\breve z_i)_k})p(t+s,\widehat{(\breve z_i)_k})]d\widehat{(\breve z_i)_k}d\hat y_k\\
&&+\int_{\bar x_k}^A\int_{\hat y_i}^A[b_i(t+s,\hat z_i)p(t+s,\hat z_i)]d\hat z_idy\Big)ds.
  \eess
Using (\ref{3.19}), we obtain
  \bess
|(N_1(p))_{x_k}|&\leq&\varepsilon Q+(A_1\cdots A_{k-1}A_{k+1}\cdots A_d)^2A_kQ\\
&&+(G^3T+1)\Big[GT(A_1\cdots A_{k-1}A_{k+1}\cdots A_d)^2A_kQ\\
&&+\sum_{i=1,i\neq k}^d(A_1\cdots A_{i-1}A_{i+1}\cdots A_{k-1}A_{k+1}\cdots A_d)^2A_kSQ\\
&&+(A_1\cdots A_{i-1}A_{i+1}\cdots A_d)^2SQ\\&&+
\sum_{i,j=1,i\neq j\neq k}^dA_iA_jA_k(A_1\cdots A_{i-1}A_{i+1}\cdots A_{j-1}A_{j+1}\cdots A_{k-1}A_{k+1}\cdots A_d)^2SQ\\
&&+2\sum_{j,k=1,k\neq j}^dA_kA_j(A_1\cdots A_{j-1}A_{j+1}\cdots A_{k-1}A_{k+1}\cdots A_d)^2SQ\\
&&+
\sum_{i=1,i\neq k}^dA_iA_k(A_1\cdots A_{i-1}A_{i+1}\cdots A_{k-1}A_{k+1}\cdots A_d)^2SQ\\
&&+(A_1\cdots A_{i-1}A_{i+1}\cdots A_d)^2SQ\Big]\\
&\leq&\varepsilon Q+(1-\varepsilon)Q, \ \ \forall t\in[0,T],\ x\in\mathbb{R}^d.
  \eess
Consequently, $N_1:\ D_1\to D_1$. It follows from the above estimates that if $p_n\to p$ in sense of the topology
of the set $D_1$ we have $N_1(p_n)\to N_1(p)$ in sense of the topology
of the set $D_1$. Therefore the operator $N_1:\ D_1\to D_1$ is a continuous operator.
It follows from the definitions of $D_1$ and $\tilde D_1$ that
$D_1$ is a compact set in the space $\tilde D_1$. $\Box$

{\bf Proof of Theorem \ref{t3.4}}
Combining the Lemmas \ref{l3.7} and \ref{l3.8}, and using Lemma \ref{l3.6}, we
deduce that the operator $L_1$ has a fixed point $p^1\in D_1$. Hence $p^1$ is a solution
to (\ref{3.20}) in $D_1$.
In order to get the global existence, we need to define the set
   \bess
B_2=\{x\in\mathbb{R}^d:\ A_i\leq x_i\leq 2A_i, \ \ i=1,2,\cdots,d\}.
  \eess
We consider the problem (\ref{3.20}) in $B_2$.
In order to do that, we consider the operator
   \bess
L_2(p)&=&p(t,x)+\int_x^A\int_y^Ap(t,z)dzdy-\frac{e^{-[a]T}}{1-e^{-[a]T}}\int_0^Te^{\int_t^{t+s}a(r)dr}\Big(a(t+s)
\int_x^A\int_y^A\nonumber\\
&& p(t+s,z)dzdy+\sum_{i=1}^d\int_{\bar x_i}^A\int_{\bar y_i}^A
(\sigma\sigma^T)_{ii}(t+s,\hat z_i)\big[p(t+s,\hat z_i)\\
&&-p^1(t+s,z_1,\cdots,z_{i-1},A_i,z_{i+1},\cdots,z_d)\\&&
+(A_i-x_i)p_{x_i}^1(t+s,z_1,\cdots,z_{i-1},A_i,z_{i+1},\cdots,z_d)\big]d\hat z_id\hat y_i\nonumber\\
&&+\sum_{i,j=1,i\neq j}^d\int_{x}^A\int_{\bar y_{ij}}^A
(\sigma\sigma^T)_{ij}(t+s,\hat z_{ij})p(t+s,\hat z_{ij})d\hat z_{ij}dy\nonumber\\
&&+\sum_{i=1}^d\int_x^A\int_{\breve y_i}^A[b_i(t+s,\breve z_i)p(t+s,\breve z_i)]d\breve z_idy\Big)ds
  \eess
under the sets
      \bess
D_2&=&\{u\in\mathcal{C}^1([0,T],\mathcal{C}^2(B_2)):\ |u|\leq Q,|u_t|\leq Q,|u_{x_i}|\leq Q, \ i=1,2,\cdots,d\},\\
\tilde D_2&=&\{u\in\mathcal{C}^1([0,T],\mathcal{C}^2(B_2)):\ |u|\leq(1+\varepsilon) Q,|u_t|\leq (1+\varepsilon)Q,\\
&&\qquad|u_{x_i}|\leq (1+\varepsilon)Q, \ i=1,2,\cdots,d\}.
  \eess
In the set $D_2$ and $\tilde D_2$, we define a norm as follows:
  \bess
\|u\|=\max\left\{\max_{t\in[0,T],x\in B_2}|u|,\ \max_{t\in[0,T],x\in B_2}|u_t|,\ \max_{t\in[0,T],x\in B_2}|u_{x_i}|,\ i=1,2,\cdots,d\right\}.
  \eess
Similar to the operator $L_1$, we define
   \bess
L_2(p)=M_2(p)+N_2(p),
  \eess
where
  \bess
M_2(p)&=&(1+\varepsilon)p,\\
N_2(p)&=&-\varepsilon p+\int_x^A\int_y^Ap(t,z)dzdy-\frac{e^{-[a]T}}{1-e^{-[a]T}}\int_0^Te^{\int_t^{t+s}a(r)dr}\nonumber\\
&& \Big(a(t+s)
\int_x^A\int_y^Ap(t+s,z)dzdy+\sum_{i=1}^d\int_{\bar x_i}^A\int_{\bar y_i}^A
(\sigma\sigma^T)_{ii}(t+s,\hat z_i)\big[p(t+s,\hat z_i)\\
&&-p^1(t+s,z_1,\cdots,z_{i-1},A_i,z_{i+1},\cdots,z_d)\\&&
+(A_i-x_i)p_{x_i}^1(t+s,z_1,\cdots,z_{i-1},A_i,z_{i+1},\cdots,z_d)\big]d\hat z_id\hat y_i\nonumber\\
&&+\sum_{i,j=1,i\neq j}^d\int_{x}^A\int_{\bar y_{ij}}^A
(\sigma\sigma^T)_{ij}(t+s,\hat z_{ij})p(t+s,\hat z_{ij})d\hat z_{ij}dy\nonumber\\
&&+\sum_{i=1}^d\int_x^A\int_{\breve y_i}^A[b_i(t+s,\breve z_i)p(t+s,\breve z_i)]d\breve z_idy\Big)ds.
  \eess
Similar to the case of $L_1$, we obtain the operator $L_2$ has a fixed point
$p^2(t,x)$ in the set $D_2$, which is a solution to (\ref{3.20}) in $D_2$.
Following Lemma \ref{l3.5}, we have $p^2(t,x)$ satisfies
   \bes
&&\int_x^A\int_y^Ap^2(t,z)dzdy-\frac{e^{-[a]T}}{1-e^{-[a]T}}\int_0^Te^{\int_t^{t+s}a(r)dr}\nonumber\\
&& \Big(a(t+s)
\int_x^A\int_y^Ap^2(t+s,z)dzdy+\sum_{i=1}^d\int_{\bar x_i}^A\int_{\bar y_i}^A
(\sigma\sigma^T)_{ii}(t+s,\hat z_i)\big[p^2(t+s,\hat z_i)\nonumber\\
&&-p^1(t+s,z_1,\cdots,z_{i-1},A_i,z_{i+1},\cdots,z_d)\nonumber\\&&
+(A_i-x_i)p_{x_i}^1(t+s,z_1,\cdots,z_{i-1},A_i,z_{i+1},\cdots,z_d)\big]d\hat z_id\hat y_i\nonumber\\
&&+\sum_{i,j=1,i\neq j}^d\int_{x}^A\int_{\bar y_{ij}}^A
(\sigma\sigma^T)_{ij}(t+s,\hat z_{ij})p^2(t+s,\hat z_{ij})d\hat z_{ij}dy\nonumber\\
&&+\sum_{i=1}^d\int_x^A\int_{\breve y_i}^A[b_i(t+s,\breve z_i)p^2(t+s,\breve z_i)]d\breve z_idy\Big)ds=0.
  \lbl{3.22}\ees
When $x_1=A_1$, the above equality deduces that
   \bess
&&\frac{e^{-[a]T}}{1-e^{-[a]T}}\int_0^Te^{\int_t^{t+s}a(r)dr}\Big(\int_{\bar x_1}^A\int_{\bar y_1}^A
(\sigma\sigma^T)_{11}(t+s,A_1,z_{2},\cdots,z_d)\\
&&\big[p^2(t+s,A_1,z_{2},\cdots,z_d)-p^1(t+s,A_1,z_{2},\cdots,z_d)\big]d\hat z_1d\hat y_1=0.
  \eess
Differentiating the above equality with respect to $t$, and using the periodicity of $\sigma,p^1$ and $p^2$
and $(\sigma\sigma^T)_{ij}>0$,
one can obtain
   \bess
p^1(t,A_1,x_2,\cdots,x_d)=p^2(t,A_1,x_2,\cdots,x_d).
  \eess
Differentiating (\ref{3.22}) with respect to $x_1$, after which we put $x_1=A_1$, we get
    \bess
&&\frac{e^{-[a]T}}{1-e^{-[a]T}}\int_0^Te^{\int_t^{t+s}a(r)dr}\Big(\int_{\bar x_1}^A\int_{\bar y_1}^A
(\sigma\sigma^T)_{11}(t+s,A_1,z_{2},\cdots,z_d)\\
&&\big[p_{x_1}^2(t+s,A_1,z_{2},\cdots,z_d)-p_{x_1}^1(t+s,A_1,z_{2},\cdots,z_d)\big]d\hat z_1d\hat y_1=0.
  \eess
Similarly, we get
      \bess
p_{x_1}^1(t,A_1,x_2,\cdots,x_d)=p_{x_1}^2(t,A_1,x_2,\cdots,x_d).
  \eess
As in above discussion, we can deduce that
   \bess
p^1(t,x_1,A_2,\cdots,x_d)&=&p^2(t,x_1,A_2,\cdots,x_d),\\
p_{x_2}^1(t,x_1,A_2,\cdots,x_d)&=&p_{x_2}^2(t,x_1,A_2,\cdots,x_d),\\
&\vdots&\\
p^1(t,x_1,x_2,\cdots,A_d)&=&p^2(t,x_1,x_2,\cdots,A_d),\\
p_{x_d}^1(t,x_1,x_2,\cdots,A_d)&=&p_{x_d}^2(t,x_1,x_2,\cdots,A_d).
  \eess
The function
   \bess
p(t,x)=\left\{\begin{array}{lll}
p^1(t,x), \ \ \ t\geq0,\ x\in B_1,\\
p^2(t,x), \ \ \ t\geq0,\ x\in B_2,
  \end{array}\right.\eess
is a solution to the following problem
   \bess\left\{\begin{array}{llll}
\partial_tp=-\sum_i\frac{\partial}{\partial x_i}(b_ip)+\frac{1}{2}\sum_{i,j}
\frac{\partial^2}{\partial x_i\partial x_j}((\sigma\sigma^T)_{ij}p), \ \ \ & t>0,\ \ x\in B_1\cup B_2,\\
p(t,x)=p(t+T,x), \ \ \ & t\geq0,\ \ x\in B_1\cup B_2.
  \end{array}\right.\eess
Repeating the above steps, using partitioning of
$\mathbb{R}^d$ into cubes, we obtain a periodic solution to  (\ref{3.17}).
The proof is complete. $\Box$

\medskip

\noindent {\bf Acknowledgment} The first author was supported in part
by NSFC of China grants 11771123, 11726628 and 11531006.

 \end{document}